\documentclass[12pt,nopreprintline]{elsarticle}

\usepackage{framed,multirow}
\usepackage[mathlines]{lineno}
\usepackage{amssymb,amsthm}
\usepackage{latexsym}
\usepackage{bbm}
\usepackage{ulem}
\usepackage{dsfont}
\usepackage{hyperref}
\usepackage{amsfonts}
\usepackage{graphicx}
\usepackage{subfig}
\usepackage{multicol}
\usepackage{epstopdf}
\usepackage{overpic}
\usepackage{amsmath,esint}
\usepackage{epsfig}
\usepackage{tikz}
\usepackage{cleveref}
\usetikzlibrary{arrows}
\usetikzlibrary{fadings}
\usepackage{pgfplots}
\pgfplotsset{compat=1.10}
\usepgfplotslibrary{fillbetween}
\usetikzlibrary{patterns}
\usepackage{appendix}
\usepackage{natbib}
\usepackage[utf8]{inputenc}
\usepackage{bbm}
\usepackage{enumitem}
\usepackage{mathrsfs}
\usepackage{verbatim}
\usepackage{graphicx}
\usepackage{algorithm}
\usepackage{tcolorbox}
\usepackage{algorithmicx}
\usepackage{algpseudocode}
\usepackage{mathtools}
\usepackage{pifont}
\usepackage{calc}
\usepackage{accents}
\usepackage{cancel}
\usepackage[margin=1in]{geometry}
\usepackage{ifpdf}
\usepackage{array}

\setenumerate{label=(\roman*),itemsep=3pt,topsep=3pt}

\newcommand{\RN}[1]{%
  \textup{\uppercase\expandafter{\romannumeral#1}}%
}

\DeclarePairedDelimiter{\ceil}{\lceil}{\rceil}
\setenumerate{label=(\roman*),itemsep=3pt,topsep=3pt}

\newcommand{\bbR}{\mathbb R}

\newcommand{\cC}{\mathcal C}

 \newcommand{\cJ}{\mathcal J}
 
 \newcommand{\cN}{\mathcal N}

\newcommand{\cS}{\mathcal S} 
 \newcommand{\cV}{\mathcal V}

\newcommand{\Rthree}{{\bbR^3}}
\newcommand{\Stwo}{{\cS^2}}

\graphicspath{{Fig/}}

\begin{document}

\begin{frontmatter}

\title{Adjoint DSMC for Nonlinear Spatially-Homogeneous Boltzmann Equation With a General Collision Model}


\author[address1]{Yunan Yang\corref{mycorrespondingauthor}}
\ead{yunan.yang@eth-its.ethz.ch}
\cortext[mycorrespondingauthor]{Corresponding author}
\author[address2]{Denis Silantyev}
\ead{dsilanty@uccs.edu}
\author[address3]{Russel  Caflisch}
\ead{caflisch@courant.nyu.edu}
\address[address1]{Institute for Theoretical Studies, ETH Zurich, 8092 Zurich, Switzerland.}
\address[address2]{Department of Mathematics, University of Colorado, Colorado Springs, CO 80918.}
\address[address3]{Courant Institute of Mathematical Sciences, New York University, New York, NY 10012.}

 \begin{abstract} 
We derive an adjoint method for the Direct Simulation Monte Carlo (DSMC) method for the spatially homogeneous Boltzmann equation with a general collision law. This generalizes our previous results in [Caflisch, R., Silantyev, D. and Yang, Y., 2021. Journal of Computational Physics, 439, p.110404], which was restricted to the case of Maxwell molecules, for which the collision rate is constant. The main difficulty in generalizing the previous results is that a rejection sampling step is required in the DSMC algorithm in order to handle the variable collision rate. We find a new term corresponding to the so-called score function in the adjoint equation and a new adjoint Jacobian matrix capturing the dependence of the collision parameter on the velocities. The new formula works for a much more general class of collision models.
 \end{abstract}



\end{frontmatter}



\section{Introduction} \label{sec:Intro}
Kinetic equations have gained great popularity in the past three decades as modeling tools beyond their classical application regime of statistical physics~\cite{cordier2005kinetic,pareschi2013interacting,burini2016collective}. Evolution-type equations for the statistical function of the car, human, or animal positions in a framework similar to the kinetic theory of gases have been widely used to model the dynamics of traffic, human crowds, and swarms~\cite{albi2019vehicular}. In particular, the Boltzmann equation, mostly known to model rarefied gas, has been studied as the energy-transport model for semiconductors~\cite{poupaud1991diffusion,ben1996energy}, financial Brownian motion~\cite{kanazawa2018derivation}, wealth distribution~\cite{pareschi2014wealth}, and sea ice dynamics~\cite{davis2020super}.

Due to the rise in computational power and the capability to collect an enormous amount of data, data-driven modeling has become a practical mainstream approach. An essential component of data-driven modeling is optimization, where the physical or nonphysical model parameters are regarded as the optimizer of an objective function that depends on the solution of the forward kinetic model~\cite{albi2015kinetic}. There is usually no explicit formula for the dependence between the model parameter and the solution to the Boltzmann equation~\cite{albi2014boltzmann}. Moreover, the large dimensionality of the unknown parameter makes it impractical to perform a global search to solve the corresponding optimization problem. Due to these two major challenges, one often turns to gradient-based optimization algorithms to perform large-scale model parameter calibration. A computationally efficient method to compute the gradient is thus essential, as iterative optimization algorithms such as gradient descent may require hundreds of gradient evaluations in a large-scale optimization task.

In an earlier work~\cite{caflisch2021adjoint}, we investigated the spatially homogeneous Boltzmann equation constrained optimization problems and developed Monte Carlo algorithms to compute the gradient based on the two common approaches in the context of PDE-constrained optimization: optimize-then-discretize (OTD) and discretize-then-optimize (DTO). Since the Boltzmann equation is an integro-differential equation, so is the corresponding adjoint equation. Obtaining the gradient based on the OTD approach requires numerical solutions to both the forward and adjoint integro-differential equations. Although the forward Boltzmann equation can be solved efficiently by the so-called Direct Simulation Monte Carlo (DSMC) method~\cite{bird1970direct,nanbu1980direct,babovsky1986simulation,babovsky1989convergence,pareschi2001introduction}, solving the adjoint equation using Monte Carlo-type methods requires interpolation in the three-dimensional velocity space at each time step~\cite{caflisch2021adjoint}. As a more efficient alternative, the DTO approach computes the gradient by deriving the adjoint of the Monte Carlo discretization for the forward model. This approach gave rise to the so-called adjoint DSMC method proposed in~\cite{caflisch2021adjoint}. One highlight is that the adjoint DSMC method costs even less than the forward DSMC method since it does not require further sampling by using the same sampled velocity pairs and collision parameters from the forward DSMC.
The idea of adjoint DSMC was recently extended to a spatially inhomogeneous case in combination with the density method for topology optimization problems~\cite{guanatopology}.

The framework in~\cite{caflisch2021adjoint} is based on the spatially homogeneous Boltzmann equation for the Maxwell molecules. That is, the collision kernel in the Boltzmann equation is a constant. In this work, we generalize the adjoint DSMC method to apply to more general collision models for the Boltzmann equation that have both velocity and angle dependence. The forward DSMC method typically uses the rejection sampling method (also known as the acceptance-rejection method) for collision kernels with velocity dependence to draw velocity pairs and the collision scattering angle from a general collision kernel. As the discrete adjoint of the forward DSMC algorithm, the adjoint DSMC will also be modified to reflect the rejection sampling process. Our new generalized adjoint DSMC method is based on the Monte Carlo gradient reparameterization coupled with the rejection sampling algorithms~\cite{Blei2017,mohamed2019monte}. For collision models with angle dependence, the Jacobian matrix of the post-collision velocities with respect to the pre-collision velocities is also different from the case where the scattering angle is uniformly distributed. Combining the generalization in these two directions, the resulting back-propagation rule for the adjoint particles is a slight modification compared to the case in~\cite{caflisch2021adjoint} for the Maxwell molecules.

We can obtain the generalized adjoint DSMC algorithm through two different derivation approaches with the same final result. One is based on the direct approach by differentiating the objective function directly, and the discrete adjoint variable is interpreted as an influence function, which is the derivative with respect to the conditional expectation. The second one is a Lagrangian approach where we impose the forward DSMC collision rules through 
Lagrangian multipliers (i.e., the discrete adjoint variables) and thus treat velocity particles at all times to be independent. The Karush–Kuhn–Tucker (KKT) conditions lead to the same gradient formulation and adjoint equations as the direct approach.
The value of the two derivations is that the direct approach shows the adjoint variable is the influence function, while the Lagrangian approach should be easier to generalize to other situations. 

The derivation of the adjoint equations for collision kernels that depend on the relative velocity requires a modification of the DTO approach, and the resulting method is a mixture of DTO and OTD. As discussed in~\Cref{section:Influence} and~\Cref{sec:conclusion}, the velocity derivative (an optimization step) is applied to the expectation over the rejection sampling, and then sampling (i.e., choice of the relevant random variables) for the rejection sampling (a discretization step) is performed afterward. For all other steps of the method, the DTO approach is followed since discretization is performed before differentiation. 

The rest of the paper is organized as follows. We first present some essential background in~\Cref{sec:background} where we briefly review the spatially homogeneous Boltzmann equation with a general collision kernel and the DSMC method with provable convergence for solving such Boltzmann equations. In~\Cref{sec:adjoint DSMC}, two different approaches to deriving the adjoint DSMC method are presented in~\Cref{subsec:direct} and~\Cref{subsec:Lagrangian}, respectively. We then present the adjoint Jacobian matrix calculation for a common form of collision kernel in~\Cref{subsec:discuss D}. It is followed by a discussion in~\Cref{subsec:special} on two special cases of the collision kernel where the adjoint DSMC algorithm could be further simplified. We show two numerical examples in~\Cref{sec:numerics} for the hard sphere collision models, one without dependence on the scattering angle and one with the dependence. We demonstrate that the adjoint DSMC method computes the gradient accurately up to the standard Monte Carlo error. Conclusion follows in~\Cref{sec:conclusion}.



\section{Background}\label{sec:background}
In this section, we briefly review the spatially homogeneous Boltzmann equation and the DSMC method for a general collision kernel.

\subsection{The Spatially Homogeneous Boltzmann Equation with a General Collision Kernel}
Consider the spatially homogeneous Boltzmann equation
\begin{equation}\label{eq:homoBoltz}
    \frac{\partial f}{\partial t}  =  Q(f,f).
\end{equation}
The nonlinear collision operator $Q(f,f)$, which describes the binary collisions among particles, is defined as
\begin{equation*}
    Q(f,f) = \int_{\bbR^3} \int_{\cS^2} q(v-v_1,\sigma) (f(v_1')f(v') - f(v_1) f(v)) d\sigma dv_1,
\end{equation*}
in which  $(v',v_1')$ represent the post-collisional velocities associated with the pre-collisional velocities $(v,v_1)$, $q\geq 0$, and the $\sigma$ integral is over the surface of the unit sphere $\cS^2$.

By conserving the momentum $v+v_1$ and the energy $v^2 + v_1^2$, we have
\begin{equation}\label{BoltzmannSolution}
\begin{split}
v' &= 1/2 (v+v_1) + 1/2 |v-v_1| \sigma,  \\
v_1'&=1/2 (v+v_1) - 1/2  |v-v_1| \sigma,
\end{split}
\end{equation}
where $\sigma \in \cS^2$ is a collision parameter. We will hereafter use the shorthand notation 
\begin{eqnarray*}
 (f, f_1, f', f_1') &=&  (  f(v), f(v_1), f(v'), f(v_1') ), \\
 (\hat f, \hat f_1,\hat f',  \hat f_1') &=&  \rho^{-1}(  f(v), f(v_1), f(v'), f(v_1') )
\end{eqnarray*}
for the values of $f$ in the Boltzmann collision operator $Q$, as well as the normalized densities $\hat f$ with $\rho = \int f(v) dv$ and $\int \hat f(v) dv=1$.

Physical symmetries imply that
\begin{equation}~\label{eq:GeneralKernel}
    q(v-v_1,\sigma) =  \tilde q(|v-v_1|,\theta)
\end{equation}
where $\cos \theta = \sigma\cdot \alpha$ with $\alpha = \frac{v-v_1}{|v-v_1|}$ and $\sigma = \frac{v'-v'_1}{|v'-v'_1|} = \frac{v'-v'_1}{|v-v_1|}$. In~\cite{caflisch2021adjoint}, we rewrote~\eqref{BoltzmannSolution}, in terms of operators, as
\begin{equation}\label{eq:AB_vel_General}
        \begin{pmatrix}
    v'\\
   v_1'
    \end{pmatrix}    = A(\sigma,\alpha) \begin{pmatrix}
    v\\
    v_1
    \end{pmatrix}, \quad 
   \begin{pmatrix}
    v\\
     v_1
    \end{pmatrix}     = B(\sigma,\alpha)    \begin{pmatrix}
    v'\\
    v_1'
    \end{pmatrix},
\end{equation}
where
\begin{equation} \label{eq:AB_def}
A(\sigma ,\alpha )= \dfrac{1}{2}\begin{pmatrix}
    I+\sigma  \alpha ^T & I-\sigma  \alpha ^T\\
    I-\sigma  \alpha ^T & I+\sigma  \alpha ^T
    \end{pmatrix},\,B(\sigma ,\alpha )= \dfrac{1}{2}\begin{pmatrix}
    I+\alpha  \sigma ^T & I-\alpha  \sigma ^T\\
    I-\alpha  \sigma ^T & I+\alpha  \sigma ^T
    \end{pmatrix},
\end{equation}
where $I$ is the identity matrix in $\Rthree$. Note that $B= A^T = A^{-1}$, showing the involutive nature of the collision. We also define the matrix $C\in\mathbb{R}^6$ as the derivative of the collision outcome $(v',v_1')$ with respect to the incoming velocities $(v,v_1)$, i.e.,
\begin{equation}\label{eq:C_def}
C(\sigma,\alpha) = \dfrac{\partial (v',v_1')}{\partial (v,v_1)}.
\end{equation}

Our main assumption on the collision kernel $q$ is that there is a positive constant $\Sigma$ such that
\begin{equation}~\label{eq:UpperBound}
    q(v-v_1,\sigma) \leq \Sigma
\end{equation}
for all  $v$, $v_1$, and $\sigma$. In practice, for example, for the numerical solution by the forward DSMC method, this is not a restriction because one can take $\Sigma$ to be the largest value of $q$ for the discrete set of velocity values.

Although the formulation and analysis presented here are valid for the general collision model~\eqref{eq:GeneralKernel}, the numerical results will be restricted to collision models of the form 
\begin{equation}\label{eq:VSS/VHS}
    q(v-v_1,\sigma) =  \tilde q(|v-v_1|,\theta) = C_\kappa(\theta) |v-v_1|^\beta.
\end{equation}
We refer to~[Sec.~1.4]\cite{villani2002review} for more modeling intuition regarding this type of collision models. The general model~\eqref{eq:VSS/VHS} accommodates both the variable hard sphere (VHS) model, when $C_\kappa(\theta)$ is constant, and the variable soft sphere (VSS) model~\cite[(4.10)]{gamba2017fast} with angular dependence. 
As described in~\Cref{subsec:special}, additional efficiency is achieved by using the separable dependence of $q$ on $|v-v_1|$ and $\theta$; see~\Cref{alg:VHS_DSMC_2}. We remark that $\sigma$ should be considered as a function of $\theta, v,v_1$ as $\theta = \arccos (\sigma \cdot \alpha)$ based on~\eqref{eq:VSS/VHS}.


With the assumption~\eqref{eq:UpperBound}, the Boltzmann equation~\eqref{eq:homoBoltz} can be further written as
\begin{eqnarray}
    \frac{\partial f}{\partial t}   &= & \int_{\Rthree}\int_\Stwo \left(f'f_1' - f f_1\right)  q (v-v_1,\sigma) d\sigma dv_1 \nonumber \\
    & =& \iint f'f_1' q  d\sigma dv_1 - \iint ff_1 q  d\sigma dv_1\nonumber  \\
    & = & \iint f'f_1' q  d\sigma dv_1  + \iint ff_1 (\Sigma  - q)  d\sigma dv_1  - \iint \Sigma  ff_1 d\sigma dv_1 \nonumber\\
    &=& \iint f'f_1' q  d\sigma dv_1  + f \iint f_1 (\Sigma - q)  d\sigma dv_1  - \mu f,\label{eq:VHS1}
\end{eqnarray}
where $\rho = \int f_1 dv_1$ is the density, $A = \int_{0}^{2\pi}\int_{0}^{\pi} \sin(\theta) d\theta d\varphi = 4\pi$ is the surface area of the unit sphere, and 
\begin{equation} \label{eq:mu}
\mu = A \Sigma \rho.
\end{equation}
If we multiply both sides of~\eqref{eq:VHS1} by an arbitrary test function $\phi(v)$, then divide by $\rho$ and finally integrate over $v$, we obtain 
\begin{eqnarray*}
       \frac{\partial ( \int \phi \hat f dv)}{\partial t}  &= & \rho \iiint  \phi \hat f' \hat f_1' q  d\sigma dv_1 dv + \rho \iiint \phi \hat f \hat f_1 (\Sigma  - q)  d\sigma dv_1 dv -  \mu \int  \phi \hat f dv\\
       & = &\rho \iiint  \phi' \hat f \hat f_1 q  d\sigma dv_1 dv +\rho \iiint  \phi \hat f \hat f_1 (\Sigma  - q)  d\sigma dv_1 dv -  \mu \int  \phi \hat f dv,
\end{eqnarray*}
using $dv dv_1=dv'dv_1' $ and then interchanging notation $(v,v_1)$ and $(v',v_1')$ in the first integral.
If we apply the explicit Euler time integration from $t_k$ to $t_{k+1} = t_k + \Delta t$, we obtain
\begin{eqnarray}
     \int \phi\ \hat f (v,t_{k+1}) dv \label{eq:weakform}
       &=&   (1-\Delta t \mu ) \int \phi\ \hat f(v,t_{k}) dv + \label{eq:non_collide} \\
        &&  \iiint  \phi'  \ {\Delta t \mu \Sigma^{-1}  {q}} \  \hat f(v,t_{k})  \hat f_1(v_1,t_{k})   A^{-1} d\sigma dv_1 dv + \label{eq:collide} \\
       &&   \iiint \phi\ {\Delta t  \mu \Sigma^{-1} (\Sigma -q)} \  \hat f(v,t_{k})  \hat f_1(v_1,t_{k})  A^{-1} d\sigma dv_1 dv. \label{eq:v_collide}
\end{eqnarray}
Note that for a fixed $k$, $\hat f(v,t_k) =  f(v,t_k)  / \rho$ is the probability density of $v$, and that $ \hat f(v,t_k)  \hat f_1(v_1,t_k) A^{-1}$ is the probability density over the three variables $(v,v_1,\sigma)$. The three terms~\eqref{eq:non_collide}, \eqref{eq:collide} and~\eqref{eq:v_collide} on the right-hand side of this equation represent the sampling of the collision process, using rejection sampling, as described in the DSMC algorithm in~\Cref{subsec:DSMC}.


\subsection{The DSMC Method}\label{subsec:DSMC}

\begin{algorithm}
  \caption{DSMC Algorithm Using Rejection Sampling\label{alg:VHS_DSMC}}
\begin{algorithmic}[1]
\State  Compute the initial velocity particles based on the initial condition, $\mathcal V^0 = \{v^0_1,\dots,v^0_N\}$. 
\For{$k=0$ to $M-1$}
\State Given $\mathcal V^{k}$, choose $N_c=\ceil[\big]{ \Delta t \mu N/2}$  velocity pairs $(i_\ell ,i_{\ell_1})$ uniformly without replacement. The remaining $N-2N_c$ particles do not have a virtual (or real) collision and set $v^{k+1}_{i} = v^{k}_{i}$.
\For{$\ell = 1$ to $N_c$}
\State Sample $\sigma^k_\ell$ uniformly over $\Stwo$.
\State Compute $\theta^k_\ell = \arccos ( \sigma^k_\ell \cdot  \alpha^k_\ell )$ and $q^k_\ell=q( |v^k_{i_{\ell }}-v^k_{i_{\ell_1}}|, \theta^k_\ell)$.
\State Draw a random number $\xi^k_\ell$ from the uniform distribution $\mathcal{U}([0,1])$.
\If{ $\xi^k_\ell \leq q^k_\ell /\Sigma $}
\State Perform real collision between $ v^k_{i_{\ell }}$ and $v^k_{i_{\ell_1}} $ following~\eqref{BoltzmannSolution} and obtain $({v^k_{i_{\ell }}}', {v^k_{i_{\ell_1}}}')$.
\State Set $( v^{k+1}_{i_{\ell }} , v^{k+1}_{i_{\ell_1 }} ) = ({v^k_{i_{\ell }}}', {v^k_{i_{\ell_1}}}'$)
\Else 
\State The virtual collision is not a real collision. Set $( v^{k+1}_{i_{\ell }} , v^{k+1}_{i_{\ell_1}} ) = ({v^k_{i_{\ell }}}, {v^k_{i_{\ell_1}}}$).
\EndIf
\EndFor
\EndFor
\end{algorithmic}
\end{algorithm}

In the DSMC method~\cite{bird1970direct,nanbu1980direct,babovsky1989convergence}, we consider a set of $N$ velocities evolving in discrete time due to collisions whose distribution can be described by the distribution function $f$ in~\eqref{eq:homoBoltz}. We divide time interval $[0,T]$ into $M$ number of sub-intervals of size $\Delta t=T/M$. At the $k$-th time interval, the particle velocities are represented as
\begin{equation*}
{\cal V}^k =\{v_1, \ldots , v_N\} (t_k),
\end{equation*}
and we denote the $i$-th velocity particle in ${\cal V}^k $ as $v_i^k$. The distribution function $f(v,t)$ is then discretized by the empirical distribution
\begin{equation}\label{eq:empirical}
    f(v,t_k) \approx \frac{\rho}{N} \sum_{i=1}^N \delta(v-v_i^k), \quad k = 0,\ldots,M.
\end{equation}

We define  the total number of virtual collision pairs $N_c =\ceil[\big]{ \Delta t \mu N/2}$. Note that the number of particles having a virtual collision is $2N_c \approx \Delta t \mu N$. Thus, the probability of having a virtual collision is $\Delta t \mu$, and the probability of not having a virtual collision is $1-\Delta t \mu$. For each velocity $v_i^k\in {\cal V}^k $ in this algorithm,  there are three possible outcomes, whose probabilities are denoted by $h_j$ where $j=1,2,3$:
%
\begin{center}
\begin{tabular}{ll}
\textbf{Outcome} &   \textbf{Probability} \\
1. No virtual collision  & $h_1 = 1-\Delta t \mu$  \\
2. A real collision      & $h_2 = \Delta t \mu \  q_i^k /\Sigma  $ \\
3. A virtual, but not a real, collision  &  $h_3 = \Delta t \mu (1- q_i^k /\Sigma)$ \\
\end{tabular}
\end{center}
Here, $q_i^k = q(v_i^k-v_{i_1}^k,\sigma_i^k)$ where $v_{i_1}^k$ represents the virtual collision partner of $v_i^k$ and $\sigma_i^k$ is the sampled collision parameter for this pair. Note that the total probability of no real collision is $h_1+h_3= 1 - \Delta t \mu q_i^k/\Sigma$. Since the $h_j$'s depend on the collision kernel $q$, we may also view them each as function of $v-v_1$ and $\sigma$. In later sections, we will use the fact that
\begin{eqnarray*}
\partial_{v_i^k} (\log h_1) &=& 0, \nonumber  \\
\partial_{v_i^k} (\log h_2) &=& (q_i^k)^{-1}  \partial_{v_i^k} q_i^k, \\
\partial_{v_i^k} (\log h_3)&=& - (\Sigma -q_i^k)^{-1}  \partial_{v_i^k} q_i^k,   \nonumber
\end{eqnarray*}
since $\rho$, $\Delta t$ and $\Sigma$ are constants. Also note that $\partial_{v_{i_1}^k} (\log h_j)= - \partial_{v_{i}^k} (\log h_j)$.

We can decide whether a particle participates in a \textit{virtual} collision or not through uniform sampling. However, in order to determine whether a selected virtual collision pair participates in the  \textit{actual} collision or not, we need to use the rejection sampling since the probability $q_i^k$ is pair-dependent. We further remark that if a virtual velocity pair is rejected for a real collision (Outcome \#2), it is automatically accepted for a virtual but not real collision (Outcome \#3). With all the notations defined above, we present the DSMC algorithm using the rejection sampling in~\Cref{alg:VHS_DSMC}. The algorithm applies to any general collision kernel satisfying~\eqref{eq:UpperBound}. 

Note that the DSMC algorithm is a single sample of the dynamics for $N$ particles following equations \eqref{eq:weakform}-\eqref{eq:v_collide}. It is not the same as $N$ samples of a single particle because of the nonlinearity of these equations. For consistency, the derivation of the adjoint equations also will employ a single sample of the $N$-particle dynamics.



\section{The Adjoint DSMC Method for a General Kernel} \label{sec:adjoint DSMC}

Consider an optimization problem for the spatially homogeneous Boltzmann equation~\eqref{eq:homoBoltz}. 
The initial condition is
\begin{equation}
f(v,0)=f_0(v;m), \label{BoltzmannIC}
\end{equation}
where $f_0$ is the prescribed initial data depending on a parameter $m$. The goal is to find $m$ which optimizes the objective function at time $t=T$,
\begin{equation}\label{eq:OTD_obj}
J (m) = \int_{\bbR^3} \phi (v) f(v,T) dv,
\end{equation}
where $f(v,T)$ is the solution to~\eqref{eq:homoBoltz} given the initial condition~\eqref{BoltzmannIC}, and thus $f(v,T)$ depends on $m$ through the initial condition.

The adjoint DSMC method proposed in~\cite{caflisch2021adjoint} is an efficient particle-based method to compute the gradient of the objective function~\eqref{eq:OTD_obj} based on the forward DSMC scheme (\Cref{subsec:DSMC}). In this work, we generalize the setup in~\cite{caflisch2021adjoint} and extend the algorithm to a general collision kernel. This section presents two different ways to obtain the general adjoint DSMC algorithm: a direct approach and a Lagrangian multiplier approach.

\subsection{A Direct Approach} \label{subsec:direct}
In this subsection, we present a direct approach by directly differentiating an equivalent form of the objective function~\eqref{eq:OTD_obj} by rewriting it as the expectation over $N$ particles at time $t=T$ (see~\cref{eq:obj1} below).

\subsubsection{Expectations for DSMC with $N$ Particles}\label{section:Expectations}
Here we define the expectations for each step of the forward DSMC algorithm. For simplicity, we assume that the number of particles $N$ is even so that the number of particle pairs is $N/2$. If $N$ is odd, there is no real difference since the single unpaired particle does not get involved in any collisions. 
The velocities change at a discrete time in the DSMC~\Cref{alg:VHS_DSMC}. This involves the following two steps at time $t_k$ where $k=0,\ldots,M-1$.

The first step is to randomly (and uniformly) select collision pairs, $v_i^k$ and $v_{i_1}^k$, and also to randomly select collision parameters.
The expectation over this step will be denoted as $E_p^k$ (with ``p" signifying collision ``parameters"). 

The second step is to perform collisions at the correct rate, i.e., the given collision kernel, $q(v -v_i, \sigma)$, using the rejection sampling. This is performed by choosing outcome $j$ with probability $h^k_{ij}=h_j (v^k_i -v^k_{i_1}, \sigma_i^k)$ for $j=1,2,3$ for this collision pair $(v_i^k, v_{i_1}^k)$; see~\Cref{subsec:DSMC} for the definition of $\{h_j\}$. The expectation for this step will be denoted as $E_a^k$ (with ``a" signifying collision ``acceptance-rejection"). Given all velocity particles at $t_k$, the total expectation over the step from $t_k$ to $t_{k+1}$ is
\begin{equation*}
E^k = E_p^k E_a^k  .
\end{equation*}
The introduction of $E_a^k$ here is motivated by the related work on calculating Monte Carlo gradients where the samples are drawn using the rejection sampling~\cite{Blei2017,mohamed2019monte}.

Moreover, the  expectation operators $E_a^k$ and $E_p^k$ can each be factored into a series of  expectations,
\begin{subequations}
\begin{align}
E_a^k &= \Pi_{i=1}^{N/2}  E_{ai}^k, \label{factorizationEa} \\
E_p^k &= \Pi_{i=1}^{N/2}  E_{pi} ,\label{factorizationEp}
\end{align}
\end{subequations}
in which the  expectations $E_{ai}^k$ and $E_{pi}^k$ are the corresponding expectations applied to a single pair $(v_i^{k}, v_{i_1}^{k})$. Note that in these products, the order of expectations does not matter since their application to a pair $(v_i^{k}, v_{i_1}^{k})$ does not affect other pairs. Also, we will only be using the factorization of $E_a^k$ in (\ref{factorizationEa}).

\subsubsection{Passing Velocity Derivatives Through the Expectations}\label{section:Formulation}

For an arbitrary test function $\phi$, the expectation for acceptance-rejection for a single pair of velocities $(v_i^k, v_{i_1}^k)$  can be written  as
\begin{equation*}
E_{ai}^k [ \phi ] 
= \sum_{j=1}^3 h_{ij}^k \phi_j,
\end{equation*}
where the random variable $\phi = \phi_j$ with probability $h_{ij}^k$. Since $h_{ij}^k$ depends on $(v_i^{k}, v_{i_1}^{k})$, for both $j=2$ and $j=3$, the velocity derivative does not commute with $E_{ai}^k$, but does commute with $E_{ai'}^k$  for $i'\neq i$ and $i'\neq i_1$.
We can calculate the velocity derivative as follows:
\begin{eqnarray*} 
  \frac{\partial}{\partial v_i^k}  E_{ai}^k [ \phi  ] 
= \sum_{j=1}^3  h_{ij}^k \left(  \frac{\partial}{\partial v_i^k}  \phi_j + \phi_j  \frac{\partial}{\partial v_i^k} \log h_{ij}^k \right) 
= E_{ai}^k \left[   \frac{\partial}{\partial v_i^k}  \phi  + \phi  \frac{\partial}{\partial v_i^k} \log h_{i}^k \right].
\end{eqnarray*}
In the last term, $h_i^k$ can be considered as a random variable for the velocity pair $(v_i^{k}, v_{i_1}^{k})$, and $h_i^k = h_{ij}^k$ with probability $h_{ij}^k$. Since $E^k_a= \Pi_{i'=1}^{N/2} E^k_{ai'}$ and ${\partial}/{\partial v_i^k}$ commutes with $E^k_{ai'}$ for $i'\neq i$ and $i'\neq i_1$, it follows that 
\begin{equation*}
\frac{\partial}{\partial v_i^k}  E_{a}^k [ \phi  ] 
= E_{a}^k \left[ \frac{\partial}{\partial v_i^k}  \phi + \phi\, \left(\frac{\partial}{\partial v_i^k} \log h_i^k \right)   \right]  .
\end{equation*}

Because $E_p^k$ also commutes with velocity derivatives, 
it follows that
\begin{equation}\label{eq:phi_diff_final}
    \frac{\partial}{\partial v_i^k}  E^k [ \phi  ]   =  \frac{\partial}{\partial v_i^k}  E_{p}^k E_{a}^k [ \phi  ]  = 
    E^k \left[ \frac{\partial}{\partial v_i^k}  \phi + \phi \left(\frac{\partial}{\partial v_i^k} \log h_i^k \right) \right].
\end{equation}

\subsubsection{Objective Function}

We rewrite the objective function~\eqref{eq:OTD_obj} at time $T= t_M$ as
\begin{equation} \label{eq:obj1}
J = E\left[ \bar{\phi}^M \right],\quad \text{where} \quad 
\bar{\phi}^M = \frac{\rho}{N}\sum_{i=1}^{N} \phi_i^M,
\end{equation}
where $\rho = \int f(v,T) dv$, $\phi_i^M=\phi(v_i^M)$ with $v_i^M \sim f(v,T)$ (in which the notation  $v \sim f$ means that the random variable  $v$ is sampled from the distribution function $f$).
Note that~\eqref{eq:obj1} is an average of expectations over the $N$ velocities, rather than a single expectation of $v$ with respect to a density function $f(v, T)$ as defined in~\eqref{eq:OTD_obj}. 

Next, we denote the conditional average at time $t_k$ as 
\begin{equation} \label{eq:Jk}
J^k = E\left[ \bar{\phi}^M | {\cal V}^k\right],
\end{equation}
which is the expectation of $\bar{\phi}^M$ for known values of the velocities $v_i^k$ at time $t_k$. 
Since
\begin{equation*}
E[\,\cdot\, | {\cal V}^k] = E^k  \ldots E^{M-1}[\,\cdot\,  | {\cal V}^k], 
\end{equation*}
we then have
\begin{equation*}
E\left[ \bar{\phi}^M | {\cal V}^k \right] = E^k \left[ \ E\left[ \bar{\phi}^M | {\cal V}^{k+1} \right] \ | {\cal V}^k \right],  
\end{equation*}
which implies
\begin{equation} \label{iteratedJ}
J^k = E^k\left[ J^{k+1} | {\cal V}^k \right],\quad k =0,\ldots, M-1.
\end{equation}

If we also define 
\begin{equation}\label{eq:Jki}
   J_i^k = E\left[ \phi_i^M | {\cal V}^k\right] ,
\end{equation}
then by the same reasoning
\begin{equation*}
J^k_i = E^k\left[ J^{k+1}_i | {\cal V}^k \right],\quad k =0,\ldots, M-1.
\end{equation*}
Moreover, for a fixed index $i=1,\ldots,N$, we have the following
\begin{equation*}
 J_i^k  = E[\phi_i^M | \cV^k] = E_{pi}^k E_{ai}^k \ldots E_{pi}^{M-1} E_{ai}^{M-1} [\phi_i^M | \cV^k], \quad k =0,\ldots, M-1.
\end{equation*}
Recall that $\{E_{pi}^k\}$ are expectations over the choice of collision pairs and parameters, respectively, and $\{E_{ai}^k\}$ are the expectations over the rejection sampling at each time step for the $i$-th particle as defined in~\Cref{section:Expectations}. Each probability $h_{ij}^k$, $j=1,2,3$, at time $t_k$ for particle index $i$, depends on $v_i^k$ and its collision partner $v_{i_1}^k$ at time $t_k$. Thus, for any $k$, 
\begin{equation} \label{eq:Jk_2}
 J^k  = \frac{\rho}{N}\sum_{i=1}^N J^k_i = \frac{\rho}{N} \sum_{i=1}^N  E[\phi_i^M | \cV^k] = \frac{\rho}{N} \sum_{i=1}^N  E_{pi}^k E_{ai}^k \ldots E_{pi}^{M-1} E_{ai}^{M-1} [\phi_i^M | \cV^k],
\end{equation}
which is another way to represent~\eqref{eq:Jk}.

Note that in~\eqref{factorizationEa}-\eqref{factorizationEp}, the index $i$ ranges from $1$ to $N/2$ as we treat $(i,i_1)$ as a pair in the rejection sampling process. In~\eqref{eq:Jk_2}, the index for the particle ranges from $1$ to $N$ because we need to treat the $i$-th and the $i_1$-th particles separately, since the objective function~\eqref{eq:obj1} is a sum of $N$ terms. Thus, in~\eqref{eq:Jk_2}, we should understand $E_{pi_1}^k E_{ai_1}^k = E_{pi}^k E_{ai}^k$ if $(i,i_1)$ is a collision pair at time $t_k$ where $i$ ranges from $1$ to $N/2$ as in~\eqref{factorizationEa}-\eqref{factorizationEp}.

\subsubsection{Influence Function}\label{section:Influence}
We introduce the influence function $\gamma_i^k=\gamma_i^k({\cal V}^k)$ defined as
\begin{equation}\label{eq:gamma1}
\gamma_i^k := \frac{\partial}{\partial v_i^k}J^k = \frac{\partial}{\partial v_i^k}  E^k [ J^{k+1} | {\cal V}^k ] =  \frac{\rho}{N }\sum_{i'=1}^N \frac{\partial  }{\partial v_i^k}  E^k[   J_{i'}^{k+1} | {\cal V}^k ],
\end{equation}
where the second equality is a result of~\eqref{iteratedJ}, and third equation follows \eqref{eq:Jk_2}. 
After applying~\eqref{eq:phi_diff_final}, we have
\begin{equation}\label{eq:gamma2}
\frac{\partial  }{\partial v_i^k} J_{i'}^{k}  = \frac{\partial  }{\partial v_i^k}  E^k [ J_{i'}^{k+1} | {\cal V}^k ] = E^k \left[ \frac{\partial}{\partial v_i^k}  J_{i'}^{k+1}  + J_{i'}^{k+1} \left(\frac{\partial}{\partial v_i^k} \log h_{i'}^k \right) \Big| {\cal V}^k   \right],
\end{equation}
for $i' = 1,\ldots, N$. Note that $h_{i'j}^k = h_{i_1' j}^k$ if $(i',i_1')$ is a pair at $t_k$, for $j=1,2,3$.
Based on the forward DSMC algorithm, the acceptance-rejection probabilities $h_{i'}^k$ in $\{J_{i'}^k\}$  do not depend on $v_i^k$ if $i' \neq i, i_1$. Thus, the last term in~\eqref{eq:gamma2} disappears for $i' \neq i, i_1$. 
After calculating the sum in~\eqref{eq:gamma1}, we have
\begin{equation}
\gamma_i^k =  E^k\left[  \frac{\partial}{\partial v_i^k } J^{k+1}  + \frac{\rho}{N}J^{k+1}_i \left(\frac{\partial}{\partial v_i^k}\log h^k_i \right)  + \frac{\rho}{N} J^{k+1}_{i_1} \left(\frac{\partial}{\partial v_i^k}\log h^k_i \right) \Big| {\cal V}^k  \right], \label{iteratedGamma}
\end{equation}
where $(i,i_1)$ is a pair at $t=t_k$ in the forward DSMC.



Note that $J^k$ is a function of ${\cal V}^k = \{v_i^k\}$ and  $J^{k+1}$ is a function of ${\cal V}^{k+1} = \{v_i^{k+1}\}$. 
The only elements in ${\cal V}^{k+1}$ that depend on $v_i^k$ are $v_i^{k+1}$ and $v_{i_1}^{k+1}$, which are determined by 
\begin{equation} \label{eq:general_collision}
\begin{pmatrix}
v_i^{k+1}\\
v_{i_1}^{k+1}
\end{pmatrix} = 
\mathcal{C}^k_{i} 
\begin{pmatrix}
v_i^{k}\\
v_{i_1}^{k}
\end{pmatrix},
\end{equation}
where the operator $\mathcal{C}_{i}^k$ is one of the following three possible matrices each with probability $h_{ij}^k$, $j=1,2,3$; see~\eqref{eq:AB_def} for the definition of $A(\sigma,\alpha)$ and~\Cref{section:Formulation} for details about $\{h_{ij}^k\}$.
\begin{equation}\label{eq:collision_operator}
{\cal C}_i^k = 
    \begin{cases}
    I,  & j=1, \text{$(v_i^k,v_{i_1}^k)$ does not have a real or virtual collision}, \\
     A(\sigma_i^k, \alpha_i^k), & j=2,  \text{$(v_i^k,v_{i_1}^k)$ has a real collision},\\
    I, & j=3,  \text{$(v_i^k,v_{i_1}^k)$ has a virtual, but not a real, collision},
    \end{cases}
\end{equation}
where $I \in \mathbb{R}^{6\time 6}$ is the identity matrix. Note that, in~\eqref{eq:general_collision}, $v_i^{k}$ and $v_{i_1}^{k}$ are a collision pair, but that $v_i^{k+1}$ and $v_{i_1}^{k+1}$ are not a collision pair. 

Furthermore, based on the definition of $C$ in~\eqref{eq:C_def}, we have
\begin{equation} \label{eq:back_D}
    D =  \left[ \dfrac{\partial (v',v'_{1})}{\partial (v,v_1)} \right]^\top = \begin{cases}
    I,  & j=1, \\
    [C(\sigma, \alpha)]^\top , & j=2,\\
    I, & j=3.
    \end{cases}
\end{equation}
We remark here that $A(\sigma, \alpha)$ and $C(\sigma, \alpha)$ may not be the same. We will discuss the concrete form of $D$ in detail in~\Cref{subsec:discuss D}.
It follows that 
\begin{equation}\label{backPropagation}
 \begin{pmatrix}
\partial/\partial v_i^{k}\\
\partial/\partial v_{i_1}^{k}
\end{pmatrix}
=  \frac{\partial (v_i^{k+1},v_{i_1}^{k+1})}{\partial (v_i^{k},v_{i_1}^{k})} \cdot 
 \begin{pmatrix}
\partial/\partial v_i^{k+1}\\
\partial/\partial v_{i_1}^{k+1}
\end{pmatrix}
 =  D_i^k  
 \begin{pmatrix}
\partial/\partial v_i^{k+1}\\
\partial/\partial v_{i_1}^{k+1}
\end{pmatrix},
\end{equation}
where $D_i^k$ is the matrix $D$ defined in~\eqref{eq:back_D} applied to the pair $(v_i^k, v_{i_1}^k)$.

Now we can calculate the change in the influence function over a single step from $t_k$ to $t_{k+1}$. Similar to~\eqref{iteratedGamma}, we get $\gamma_{i_1}^k$ where
\begin{equation}
\gamma_{i_1}^k =  E^k \left[ \frac{\partial}{\partial v_{i_1}^k}  J^{k+1} + \frac{\rho}{N} J^{k+1}_{i_1} \left(\frac{\partial}{\partial v_{i_1}^k}\log h^k_i \right)  + \frac{\rho}{N} J^{k+1}_{i} \left(\frac{\partial}{\partial v_{i_1}^k}\log h^k_i \right) \Big| {\cal V}^k  \right]. \label{iteratedGamma1}
\end{equation}
Combining~\eqref{iteratedGamma} and \eqref{iteratedGamma1}, and applying \eqref{backPropagation}, we have
\begin{equation}\label{eq:direct_adjoint_exp}
    \begin{pmatrix}
    \gamma_i^k\\
    \gamma_{i_1}^k
    \end{pmatrix}
     = E^k\left[ D_i^k  \begin{pmatrix}
    \gamma_i^{k+1}\\
    \gamma_{i_1}^{k+1}
    \end{pmatrix}  \Big| {\cal V}^k \right]   + \frac{\rho}{N} E^k\left[ \left( J_{i}^{k+1} + J_{i_1}^{k+1} \right) \begin{pmatrix}
    \partial \log h_i^k /\partial v_i^k\\
    \partial \log h_i^k /\partial v_{i_1}^k\
    \end{pmatrix}  \Big| {\cal V}^k \right].
\end{equation}
Again, $h_i^k$ in the right-hand side of~\eqref{eq:direct_adjoint_exp} should be seen as a random variable, and $h_i^k = h_{ij}^k =  h_j$ with probability $h_j$ where $\{h_j\}$ is defined in~\Cref{subsec:DSMC}. 



Now we apply sampling to this expectation~\eqref{eq:direct_adjoint_exp}, as in the forward DSMC algorithm, over a single sample of the dynamics of $N$ discrete particles. We remark that based on~\eqref{eq:Jki}, $J_{i}^{k+1}$ and $J_{i_1}^{k+1}$ become $\phi_{i}^{M}$ and $\phi_{i_1}^{M}$ after sampling. Finally, we obtain 
\begin{equation}\label{eq:direct_adjoint_sample}
    \begin{pmatrix}
    \gamma_i^k\\
    \gamma_{i_1}^k
    \end{pmatrix}
     =  D_i^k  \begin{pmatrix}
    \gamma_i^{k+1}\\
    \gamma_{i_1}^{k+1}
    \end{pmatrix}   +  \frac{\rho}{N} \left( \phi_{i}^{M} + \phi_{i_1}^{M} \right)  \frac{\partial \log h_i^k }{\partial v_i^k}\begin{pmatrix}
    1\\
    -1
    \end{pmatrix},
\end{equation}
because $  \frac{\partial }{\partial v_i^k} \log h_i^k  = - \frac{\partial }{\partial v_{i_1}^k} \log h_{i}^k $ due to the symmetry of elastic binary collision. 

Note that the velocity derivative (an optimization step) is applied to the expectation in \eqref{eq:Jk_2} and then sampling (a discretization step) is applied to the result in \eqref{eq:direct_adjoint_sample}. As discussed in~\Cref{sec:Intro}, this part of the derivation follows the OTD approach. On the other hand, instead of using a probability density function $f$ as in~\eqref{eq:OTD_obj}, the expectation in~\eqref{eq:Jk_2} is based on $N$ discrete velocity particles (using particle discretization). Thus, this part follows the DTO approach. The derivation here is a combination of both DTO and OTD approaches.


\subsubsection{Final Data }
The ``final data" for $\gamma_i^k$ with $k=M$ is simply 
\begin{equation}\label{eq:direct_adjoint_final}
\gamma_i^M =  \frac{\partial}{\partial v_i^M} J^M  = \frac{\partial}{\partial v_i^M} E[\bar{\phi}^M  |{\cal V}^M] =   \frac{\partial}{\partial v_i^M} \bar{\phi}^M = \frac{\rho}{N} \phi'(v_i^M).
\end{equation}

\subsubsection{Resulting Adjoint System} 
Combining both~\eqref{eq:direct_adjoint_sample} and~\eqref{eq:direct_adjoint_final}, the resulting system is
\begin{subequations}\label{eq:gamma_rules}
\begin{align}
\gamma_i^M &= \frac{\rho}{N} \phi'(v_i^M) \label{eq:gamma_final_disc} \\
\begin{pmatrix} {\gamma_{i}^{k}}\\  {\gamma_{i_1}^{k}}\end{pmatrix} &=
D_i^k  \begin{pmatrix}
{\gamma_{i}^{k+1}}\\  {\gamma_{i_1}^{k+1}}  \end{pmatrix} +  \begin{pmatrix} {\eta_{i}^{k}}\\  {\eta_{i_1}^{k}}\end{pmatrix} \quad \mbox{for $k = 0,\ldots, M-1$} ~\label{eq:gamma_non_final_disc} 
\end{align}
\end{subequations}
in which $\eta_{i}^k$ is defined as
\begin{equation}\label{eq:back_eta}
\eta_i^k = 
    \begin{cases}
    0,  & \text{if $v_i^k$ does not have a virtual or real collision}, \\ 
    \quad \\
    \dfrac{\rho}{N} \dfrac{\partial \log q_i^k }{ \partial {v_i^k}} ( \phi^M_i + \phi^M_{i_1}), &  \text{if $v_i^k$ has a virtual collision that is a real collision},\\
     \quad \\
      \dfrac{\rho}{N}  \dfrac{\partial \log (\Sigma - q_i^k) }{ \partial {v_i^k}} ( \phi^M_i + \phi^M_{i_1}), &   \text{if $v_i^k$ has a virtual, but not a real, collision}.
    \end{cases}
\end{equation}
This can be used to calculate sensitivities.  Assuming the $L^2$ inner product for $m$, we have
\begin{equation*}
\nabla_m J = \nabla_m  \left( E_{\cV^0 }\left[ E\left[\bar{\phi}^M | \cV^0 \right] \right] \right) \approx  \sum_{i=1}^N \frac{\partial v^0_i}{\partial m}  \cdot \left(\frac{\partial}{\partial v_i^0}   E\left[\bar{\phi}^M | \cV^0 \right] \right)
\approx  \sum_{i=1}^N  \frac{\partial v_{i}^0}{\partial m}\cdot \gamma_i^0.
\end{equation*}
We denote the gradient calculated through the adjoint approach as  $\nabla_m J^{AD}$ where
\begin{equation}\label{eq:adjoint J^{AD}}
\nabla_m J^{AD} =  \sum_{i=1}^N \nabla_m  v_{i}^0  \cdot \gamma_i^0.
\end{equation}
Note that $E_{\cV^0 }$ denotes the expectation over all elements in ${\cV^0 }$, and $\forall v_i^0\in \cV^0 $, $ v_i^0\sim f_0(v;m)$, the initial distribution introduced in~\eqref{BoltzmannIC}. The approximation in~\eqref{eq:adjoint J^{AD}} corresponds to the so-called pathwise gradient estimator~\cite[Section 5]{mohamed2019monte}.

To compute the gradient represented by the last term in~\eqref{eq:adjoint J^{AD}}, we need to ``back-propagate'' \eqref{eq:gamma_rules} from $t=t_M = T$ all the way to $t=t_0 = 0$. We remark that the only differences between~\eqref{eq:gamma_rules} and the adjoint DSMC algorithm designed for the Maxwell molecules~\cite{caflisch2021adjoint} are the ``$\eta$'' terms in~\eqref{eq:gamma_non_final_disc}. We summarize the new adjoint DSMC algorithm in~\Cref{alg:adjoint_VHS}.

\begin{algorithm}[ht!]
\caption{Solve the discrete adjoint DSMC system and compute the gradient~\eqref{eq:adjoint J^{AD}}.~\label{alg:adjoint_VHS}}
\begin{algorithmic}[1]
\State Given the final-time velocities $\cV_{M}$ from the forward DSMC, set $\gamma_i^{M} = \frac{\rho}{N}\phi'(v^{M}_i)$ for  all $i$.
\For{$k=M-1$ to $0$}
\State Given $\{\gamma_{1}^{k+1},\ldots, \gamma_{N}^{k+1}\}$ from the previous iteration, the collision history and collision parameters from the forward DSMC process.
\If {$v_i^k\in {\cal V}_{k}$ did not virtually collide at $t_k$}
\State Set $\gamma_{i}^k = \gamma_i^{k+1}$.
\ElsIf {$v_i^k, v_{i_1}^k\in {\cal V}_{k}$ virtually collided at $t_k$}
\State Set $\begin{pmatrix}  \gamma_i^k \\  \gamma_{i_1}^k     \end{pmatrix}  =   D_i^k \begin{pmatrix}   \gamma_{i}^{k+1} \\  \gamma_{i_1}^{k+1}  \end{pmatrix} + {\begin{pmatrix}  \eta_{i}^k \\  \eta_{i_1}^k     \end{pmatrix}}$, where $\eta_{i}^k,\eta_{i_1}^k$ follow~\eqref{eq:back_eta}.
\EndIf
\EndFor
\State Compute the gradient following~\eqref{eq:adjoint J^{AD}}.
\end{algorithmic}
\end{algorithm}

\subsection{A Lagrangian Approach}\label{subsec:Lagrangian}
In addition to the direct derivation of the adjoint equations described above, it is useful to have a ``Langrangian" form from which the forward DSMC and the adjoint equations can be derived. 
Using the earlier notation $E^k$ in~\Cref{subsec:direct} and the objective function $J$ defined in~\eqref{eq:obj1}, we can write the Lagrangian $\mathcal J$ as 
\begin{eqnarray}
\cJ &= &  J + \frac{1}{2}\sum_{k=0}^{M-1} \sum_{i=1}^N E^k \left[   \begin{pmatrix} \gamma_i^{k+1} \\ \gamma_{i_1}^{k+1} \end{pmatrix} \cdot  \left( \cC_{i}^k \begin{pmatrix} v_i^k \\ v_{i_1}^k \end{pmatrix} - \begin{pmatrix} v_i^{k+1} \\ v_{i_1}^{k+1} \end{pmatrix} \right)   \bigg| \begin{pmatrix} v_i^k \\ v_{i_1}^k \end{pmatrix} \right] +  \nonumber \\
&&   \sum_{i=1}^N E_{\hat{v}_{0i} (m) \sim f_0(v;m)} \left[  \gamma_{i}^0 \cdot  \left( \hat{v}_{0i} (m) - v_i^0 \right) \right]. \label{eq:full_Lag}
\end{eqnarray}
The ``$1/2$'' scaling is to avoid enforcing the collision rule twice. Again, note that $v_i^{k}$ and $v_{i_1}^{k}$ are a collision pair, but that neither $(v_i^{k+1}, v_{i_1}^{k+1})$ nor $(v_i^{M}, v_{i_1}^{M})$  are a collision pair.

In contrast to the direct approach in which the velocities $\{v_i^k\}$ are chosen according to the forward DSMC and the $\{\gamma_i^{k}\}$ are chosen by the adjoint equations, in the Lagrangian approach described here, the velocities $\{v_i^k\}$ and the Lagrangian multipliers $\{\gamma_i^{k}\}$ are considered to be any random variables. We then derive the forward and adjoint DSMC equations by requiring that $\mathcal J$ is stationary with respect to variations in the  $v_i^k$'s and the $\gamma_i^{k}$'s.

\subsubsection{Collision Rules and Initial Data}
The collision rules are derived from  the derivatives of $\mathcal J$  with   respect to $ \gamma_i^{k+1}$ and $\gamma_{i_1}^{k+1}$ for $k = 0,\ldots, M-1$. Setting these to zero implies that 
\begin{equation}  \label{eq:derived_collision}
\cC_{i}^k  \begin{pmatrix} v_i^k \\ v_{i_1}^k \end{pmatrix} = \begin{pmatrix} v_i^{k+1} \\ v_{i_1}^{k+1} \end{pmatrix} .
\end{equation}
See~\eqref{eq:collision_operator} for the three possible cases of $\cC_{i}^k$,
which are the DSMC equations for a non-virtual collision, a real collision and a virtual, but not real, collision, respectively.

Similarly,  setting to zero  the derivatives  of $\mathcal J$  with   respect to $ \gamma_i^{0}$ implies that 
 \begin{equation*} 
v_i^0  = \hat{v}_{0i}(m),
\end{equation*}
which is the initial data for the forward DSMC. We assume that $\hat{v}_{0i}(m) \sim f_0(v;m)$.

\subsubsection{Parameter $m$}
If we take the derivative of $\mathcal J$ with respect to the parameter $m$, we obtain the gradient based on sampled $\{\hat{v}_{0i}\}$:
\begin{equation}\label{eq:J_m}
 \frac{\partial J}{\partial m}  = \frac{\partial \mathcal{J}}{\partial m}  \approx  \sum_{i=1}^N   \gamma_{i}^0 \cdot  \frac{\partial \hat{v}_{0i}(m) }{\partial m}.
\end{equation}
Note that~\eqref{eq:J_m} is the same as the $ \nabla_m J^{AD}$ in~\eqref{eq:adjoint J^{AD}}.

\subsubsection{Final Data}
For each $i$,  we take the derivative of $\mathcal J$ with respect to the final velocity particle $v_i^M$. We need to consider the objective function $J$ as a function of $\cV^{M-1}$. Thus, $J$ here is replaced by $J^{M-1}$ defined in~\eqref{eq:Jk}.
\begin{eqnarray*}
\frac{\partial \mathcal{J}}{\partial v_i^M} 
=  \frac{\rho}{N}  E^{M-1} \left[   \phi' (v_i^M)  |  \mathcal{V}^{M-1} \right] -  E^{M-1}  \left[ \gamma_i^M | \mathcal{V}^{M-1}\right] =   E^{M-1}  \left[ \frac{\rho}{N}\phi' (v_i^M) - \gamma_i^M  \Big| \mathcal{V}^{M-1} \right]. 
\end{eqnarray*}
Setting this derivative to $0$ implies the final condition
$ \gamma_i^M= \frac{\rho}{N} \phi' (v_i^M) $ as seen in~\eqref{eq:gamma_final_disc}.

\subsubsection{Adjoint Equations}
The first term $J$ in the Lagrangian is defined as an expectation of $\phi^M$; see~\eqref{eq:obj1}. That is,
\begin{equation*}
J  = E[\phi^M] = \rho N^{-1} \sum^{N}_{i'=1}E^{1} \ldots E^{M-1} [\phi_{i'}^M].
\end{equation*}
Since $v^k_i$  only appears in the terms $h^k_{ij}=h_j (v^k_i -v^k_{i_1}, \sigma_i^k)$, $j=1,2,3$, in $E^k$ (see~\Cref{subsec:DSMC}), then using the commutator between the derivative and the expectation as in~\eqref{eq:phi_diff_final}, we have
\begin{eqnarray}\label{eq:Jsampled}
\partial_{v^k_i}J  &=& \rho N^{-1} \sum^{N}_{i'=1}E^{1}\ldots E^{k-1} (\partial_{v^k_i}E^{k})E^{k+1}\ldots E^{M-1} [\phi_{i'}^M] \nonumber  \\
&=& \rho N^{-1} \sum^{N}_{i'=1}E^{1}\ldots E^{M-1} \left[ (\delta_{i' i}+\delta_{i' i_1})\, \left( \partial_{v^k_i} \log h_i^k  \right)\,\phi_{i'}^M \right]  \nonumber 
\\
&=& \rho N^{-1} E^{1}\ldots E^{M-1}   \left[ \partial_{v^k_i} \log h_i^k  \left(\phi_{i}^M+\phi_{i_1}^M \right)\right] 
\nonumber 
\\
&\approx& \rho N^{-1} \left( \partial_{v^k_i} \log h_i^k  \right) [\phi_{i}^M+\phi_{i_1}^M], 
\end{eqnarray}
in which the  approximation in the last step is from sampling. The $h_{i}^k$ term in all  equations above except the last one should be considered as a random variable and $h_i^k = h_{ij}^k$ with probability $h_{ij}^k$, $j=1,2,3$. We use the final term in~\eqref{eq:Jsampled} as the contribution from $J$ to the (approximate) optimality condition for $ \mathcal{J}$.




We continue by differentiating the remaining term ${\mathcal J} - J$ in~\eqref{eq:full_Lag} 
to obtain
\begin{eqnarray}
\frac{\partial (\mathcal{J}-J)}{\partial v_i^k} 
&=&   E^{k}\left[  
\begin{pmatrix} \gamma_i^{k+1} \\ \gamma_{i_1}^{k+1} \end{pmatrix}  \cdot 
\bigg\{  \cC_i^k \begin{pmatrix}
v_i^{k} \\ v_{i_1}^{k}
\end{pmatrix}  - \begin{pmatrix}
v_i^{k+1} \\ v_{i_1}^{k+1}
\end{pmatrix}  \bigg\} 
   \frac{\partial \log h_i^k}{\partial v_i^k}  \Big| \mathcal{V}^k \right] 
   +  \nonumber  \\
&&     E^{k}\left[  \overline{D}_i^k \begin{pmatrix}
\gamma_i^{k+1} \\ \gamma_{i_1}^{k+1}
\end{pmatrix} \bigg | \mathcal{V}^k \right] -    E^{k-1} \left[  \gamma_i^k  | \mathcal{V}^{k-1}  \right], \nonumber \\ 
&\approx&
\begin{pmatrix} \gamma_i^{k+1} \\ \gamma_{i_1}^{k+1} \end{pmatrix}  \cdot 
\bigg\{  \cC_i^k \begin{pmatrix}
v_i^{k} \\ v_{i_1}^{k}
\end{pmatrix}  - \begin{pmatrix}
v_i^{k+1} \\ v_{i_1}^{k+1}
\end{pmatrix}  \bigg\} 
   \frac{\partial \log h_i^k}{\partial v_i^k}    
   +      \overline{D}_i^k \begin{pmatrix}
\gamma_i^{k+1} \\ \gamma_{i_1}^{k+1}
\end{pmatrix}  -      \gamma_i^k  , \label{eq:J_v_non_final} 
\end{eqnarray}
\begin{eqnarray}
\frac{\partial (\mathcal{J}-J)}{\partial v^k_{i_1}}
&=& E^k\left[ \begin{pmatrix} \gamma_i^{k+1} \\ \gamma_{i_1}^{k+1} \end{pmatrix}  \cdot
\bigg\{  \cC_i^k  \begin{pmatrix}
v_i^{k} \\ v_{i_1}^{k}
\end{pmatrix}  - \begin{pmatrix}
v_i^{k+1} \\ v_{i_1}^{k+1}
\end{pmatrix}  \bigg\} 
   \frac{\partial \log h_i^k}{\partial v_{i_1}^k}  \bigg| \mathcal{V}^k \right] +  \nonumber  \\
&&    E^k \left[  \underline{D}_i^k \begin{pmatrix}
\gamma_i^{k+1} \\ \gamma_{i_1}^{k+1} \end{pmatrix} \bigg | \mathcal{V}^k \right] -  E^{k-1} \left[  \gamma_{i_1}^k  | \mathcal{V}^{k-1}  \right]
\nonumber \\
&\approx& 
\begin{pmatrix} \gamma_i^{k+1} \\ \gamma_{i_1}^{k+1} \end{pmatrix}  \cdot 
\bigg\{  \cC_i^k \begin{pmatrix}
v_i^{k} \\ v_{i_1}^{k}
\end{pmatrix}  - \begin{pmatrix}
v_i^{k+1} \\ v_{i_1}^{k+1}
\end{pmatrix}  \bigg\} 
   \frac{\partial \log h_i^k}{\partial v_{i_1}^k}    
   +      \underline{D}_i^k \begin{pmatrix}
\gamma_i^{k+1} \\ \gamma_{i_1}^{k+1}
\end{pmatrix}  -      \gamma_{i_1}^k. 
\label{eq:J_v1_non_final}
\end{eqnarray}
Here, $D_i^k = \begin{bmatrix}
\ \overline{D}_i^k\  \\ \ \underline{D}_i^k \
\end{bmatrix}$. That is, $\overline{D}_i^k, \underline{D}_i^k \in \mathbb{R}^{3\times 6}$ are the upper and  lower blocks of the matrix $D_i^k$ defined in~\eqref{eq:back_D}. Note that the terms in brackets ``$\{ \ \}$" are always $0$ because of~\eqref{eq:derived_collision}.

Finally, we combine~\eqref{eq:Jsampled}, \eqref{eq:J_v_non_final} and~\eqref{eq:J_v1_non_final}, to obtain approximations for $\frac{\partial \mathcal{J}}{\partial v^k_{i}}$ and $\frac{\partial \mathcal{J}}{\partial v^k_{i_1}}$. Setting them to zero results in equations that are identical to~\eqref{eq:gamma_rules}; i.e., the Lagrangian approach confirms the results of the direct approach.

Note that the two derivations of adjoint equations, the first directly from the DSMC equations (\Cref{subsec:direct}) and the second from a Lagrangian (\Cref{subsec:Lagrangian}), are nearly identical in terms of the details of the calculations and the origin of the score function.

\subsection{The Calculation of $D$ in~\eqref{eq:back_D}} \label{subsec:discuss D}
In the earlier adjoint DSMC derivations, we did not specify the explicit form of the matrix $D$ in~\eqref{eq:back_D}, which is denoted as $D_i^k$ when applied to a particular collision pair $(v_i^k, v_{i_1}^k)$. It is clear that the matrix is an identity matrix when the velocity particles do not have a real collision. Thus, we focus on the case $j=2$.

Based on the binary collision formula~\eqref{BoltzmannSolution}, there are two cases when we calculate the Jacobian matrix $C$ and its adjoint $D$: (1)~$\sigma$ can be seen to be independent of $(v,v_1)$ when the collision kernel $q(v-v_1, \sigma)$ is angle-independent, and (2)~$\sigma$ depends on $(v,v_1)$ otherwise.

\textbf{Case one -- angle-independent kernel}: Despite the fact that $\sigma$ is the unit vector along the post-collision relative velocity, $\sigma$ can be regarded as uniformly distributed over $\Stwo$ as long as the collision kernel $q(v-v_1, \sigma) = \tilde q(|v-v_1|, \theta) $ does not depend on the scattering angle $\theta$. Thus, we can consider $\sigma$ to be independent of the pre-collision relative velocity $\alpha$, and thus independent of $(v,v_1)$. If $(v,v_1)$ is a real collision pair, based on~\eqref{eq:AB_vel_General}, we have
\begin{equation} \label{eq:D angle-independent}
   D =   C(\sigma,\alpha)^\top =  \left[\dfrac{\partial (v',v_1')}{\partial (v,v_1)} \right]^\top = B(\sigma,\alpha),
\end{equation}
where $B$ is defined in~\eqref{eq:AB_def}. This formula was used in~\cite{caflisch2021adjoint} where a constant collision kernel was considered.

\textbf{Case two -- angle-dependent kernel}: In this case, we need to consider the dependence of $\sigma$ on the pre-collision particle velocities, $v$ and $v_1$, since the distribution of $\sigma$ is no longer uniform over the sphere $\Stwo$, and it depends on $\alpha = \frac{v-v_1}{|v-v_1|}$.
Going through the calculations for~\eqref{eq:C_def}, we have
\begin{align*}
 \partial_{v} v'  &=  \frac{1}{2}  \left( I + \sigma\  \alpha^\top + |u| \partial_{v} \sigma \right)  \nonumber, \\
 \partial_{v_1} v'  &=  \frac{1}{2}  \left( I - \sigma\  \alpha^\top + |u|  \partial_{v_1} \sigma\right)  \nonumber, \\
\partial_{v} {v_1'}  &= \frac{1}{2}  \left( I - \sigma\  \alpha^\top  - |u|  \partial_{v} \sigma\right)  \nonumber, \\
 \partial_{v_1} v_1'  &=  \frac{1}{2}  \left( I+\sigma \  \alpha^\top - |u|  \partial_{v_1} \sigma \right),
\end{align*}
where $u = v- v_1$.  Note that $\partial_{v} \sigma = \partial_{u} \sigma= - \partial_{v_1} \sigma$. Next, we compute $\partial_{v} \sigma$ and $\partial_{v_1} \sigma$. Note that we can also write the collision formula~\eqref{BoltzmannSolution}, following~\cite[Sec.~2.2]{wang2008particle}, as
\begin{align*}
    v' &= v  + 1/2\,\Delta U,\\
    v_1' &= v_1 - 1/2\,\Delta U,
\end{align*}
where
\begin{equation*}
    \Delta U = \begin{pmatrix}
               \frac{u_x u_z}{\sqrt{u_x^2 + u_y^2}} & -\frac{u_y |u|}{\sqrt{u_x^2 + u_y^2}} &
               u_x\\
                \frac{u_y u_z}{\sqrt{u_x^2 + u_y^2}} & \frac{u_x |u|}{\sqrt{u_x^2 + u_y^2}} &
               u_y\\
               -\sqrt{u_x^2 + u_y^2} &
               0 &
               u_z 
               \end{pmatrix}
               \begin{pmatrix}
               \sin \theta \cos \varphi\\
               \sin \theta \sin \varphi\\
               \cos \varphi
               \end{pmatrix} - u.
\end{equation*}
Thus, using the definition of $\sigma$, we have
\begin{equation}\label{eq:sigma}
    \sigma = \frac{v'- v_1'}{|v'-v_1'|} = \frac{v'- v_1'}{|v-v_1|} = \begin{pmatrix}
               \frac{u_x u_z}{|u|\sqrt{u_x^2 + u_y^2}} & -\frac{u_y }{\sqrt{u_x^2 + u_y^2}} &
               \frac{u_x}{|u|}\\
                \frac{u_y u_z}{|u| \sqrt{u_x^2 + u_y^2}} & \frac{u_x }{\sqrt{u_x^2 + u_y^2}} &
               \frac{u_y}{|u|}\\
               -\frac{\sqrt{u_x^2 + u_y^2}}{|u|} &
               0 &
               \frac{u_z}{|u|} 
               \end{pmatrix}
               \begin{pmatrix}
               \sin \theta \cos \varphi\\
               \sin \theta \sin \varphi\\
               \cos \varphi
               \end{pmatrix},
\end{equation}
which shows the dependence of $\sigma$ upon $v,v_1, \theta, \varphi$ explicitly. For fixed $\theta$ and $\varphi$, we have
\begin{eqnarray*}
    |u| \partial_u \sigma & = & |u| \left( \frac{\partial e_i^ \top \sigma }{\partial e_j^ \top u}\right)_{ij} = \begin{bmatrix} 
    G_1 \sigma &  G_2 \sigma &  G_3 \sigma
    \end{bmatrix}, \\
        |u|  \left( \partial_u \sigma \right)^\top &=&   \begin{bmatrix} 
   \sigma^\top G_1^\top \\  \sigma^\top G_2^\top  \\  \sigma^\top G_3^\top
    \end{bmatrix} = -  \begin{bmatrix}
       \sigma^\top G_1 \\  \sigma^\top G_2  \\  \sigma^\top G_3
    \end{bmatrix},
\end{eqnarray*}
where $\{e_i\}$ is the standard basis in $\bbR^3$, and matrices $G_j = -G_j^\top$, $j=1,2,3$, as defined below using the notation $\alpha = u/|u| = [\alpha_x, \alpha_y, \alpha_z]^\top$. 
\begin{align*}
   G_1 &= \begin{bmatrix}
 0 & \frac{\alpha_y}{\alpha_x^2 +\alpha_y^2  }&
  \frac{\alpha_x^2 \alpha_z}{\alpha_x^2 +\alpha_y^2  }
  \\
   & &\\
- \frac{\alpha_y}{\alpha_x^2 +\alpha_y^2  }&
 0 &
  \frac{\alpha_x\alpha_y\alpha_z}{\alpha_x^2 +\alpha_y^2 }
  \\
     & &\\
  -\frac{\alpha_x^2 \alpha_z}{\alpha_x^2 +\alpha_y^2} &
  -  \frac{\alpha_x\alpha_y\alpha_z}{\alpha_x^2 +\alpha_y^2 } &
   0
    \end{bmatrix} = \frac{1}{\alpha_x^2 +\alpha_y^2}\begin{bmatrix}
 0 & \alpha_y &
  \alpha_x^2 \alpha_z
  \\
   & &\\
- \alpha_y&
 0 &
  \alpha_x\alpha_y\alpha_z
  \\
     & &\\
  -\alpha_x^2 \alpha_z &
  -  \alpha_x\alpha_y\alpha_z &
   0
    \end{bmatrix},\\
    G_2 &= \begin{bmatrix}
 0 & -\frac{\alpha_x}{\alpha_x^2 +\alpha_y^2  }&
  \frac{\alpha_x\alpha_y \alpha_z}{\alpha_x^2 +\alpha_y^2  }
  \\
   & &\\
 \frac{\alpha_x}{\alpha_x^2 +\alpha_y^2  }&
 0 &
  \frac{\alpha_y^2\alpha_z}{\alpha_x^2 +\alpha_y^2 }
  \\
     & &\\
  -\frac{\alpha_x\alpha_y \alpha_z}{\alpha_x^2 +\alpha_y^2} &
  -  \frac{\alpha_y^2\alpha_z}{\alpha_x^2 +\alpha_y^2 } &
   0
    \end{bmatrix} = \frac{1}{\alpha_x^2 +\alpha_y^2}\begin{bmatrix}
 0 & -\alpha_x &
  \alpha_x \alpha_y \alpha_z
  \\
   & &\\
\alpha_x &
 0 &
  \alpha_y^2 \alpha_z
  \\
     & &\\
  -\alpha_x \alpha_y \alpha_z &
-\alpha_y^2\alpha_z &
   0
    \end{bmatrix},\\   
G_3 & = \begin{bmatrix}
 0&
 0 &
 -\alpha_x
  \\
     & &\\
 0&
 0 &
 -\alpha_y
  \\
     & &\\
  \alpha_x    &\alpha_y & 0 
    \end{bmatrix}.
\end{align*}
We define a tensor $G_{lij}$ where $G_{1ij} = G_1$, $G_{2ij} = G_2$ and $G_{3ij} = G_3$. We will then use Einstein's summation convention. The adjoint matrix $D$ in~\eqref{eq:back_D} for a real collision pair becomes
\begin{eqnarray}
    D = [C(\sigma, \alpha)]^\top &=& \frac{1}{2} 
 \begin{pmatrix}
 I + \alpha\sigma^\top  &  I - \alpha\sigma^\top \\
  I - \alpha\sigma^\top &  I + \alpha\sigma^\top
\end{pmatrix}  +
\frac{1}{2} 
\begin{pmatrix}
-(\sigma)_i G_{lij}  &  (\sigma)_i G_{lij}  \\
 (\sigma)_i G_{lij}  &  -(\sigma)_i G_{lij} 
\end{pmatrix} 
\nonumber \\
& =& B(\sigma,\alpha)  + \widetilde B(\sigma, \alpha), \label{eq:new D}
\end{eqnarray}
where $(\sigma)_i$ denotes the $i$-th element of vector $\sigma$, and
\begin{equation}\label{eq:tilde B def}
 \widetilde B(\sigma, \alpha) = \frac{1}{2} \begin{pmatrix}
-(\sigma)_i G_{lij}  &  (\sigma)_i G_{lij}  \\
 (\sigma)_i G_{lij}  &  -(\sigma)_i G_{lij} 
 \end{pmatrix}.
\end{equation}
In order to build $\widetilde B(\sigma, \alpha)$, we need $\sigma$ and $\alpha$, which are already stored in memory during the forward DSMC for the calculation of $B(\sigma, \alpha)$ in~\eqref{eq:AB_def}. Therefore, there is no additional memory requirement.

\subsection{Special Cases} \label{subsec:special}
In this section, we discuss a few special cases of the collision kernel $q(v-v_1,\sigma)$ that may result in a variation or simplification with respect to~\Cref{alg:adjoint_VHS}.

\subsubsection{The DSMC and the Adjoint DSMC Methods for a Constant Kernel}

For a constant collision kernel $q$, which is the model for Maxwell molecules, the upper bound can be taken as $\Sigma =q$. Then $\{h_j\}$ in~\Cref{subsec:DSMC} are all constant, all virtual collisions are real collisions, and all $\eta$ terms are $0$ in~\eqref{eq:back_eta}. This corresponds to the case studied in~\cite{caflisch2021adjoint}.

\subsubsection{The DSMC and the Adjoint DSMC Methods for Kernel of Type~\eqref{eq:VSS/VHS}}\label{sec:algs13_discuss}

\begin{algorithm}
 \caption{DSMC Algorithm for Collision Kernel of Form~\eqref{eq:VSS/VHS} \label{alg:VHS_DSMC_2}}
\begin{algorithmic}[1]
\State Compute the initial velocity particles based on the initial condition, $\cV^0 = \{v^0_1,\dots,v^0_N\}$. Set $N_c =\ceil[\big]{N \Delta t \mu_\kappa/2}$ where $\mu_\kappa = A_\kappa \Sigma_v \rho$ and $M=T/\Delta t$ given final time $T$. 
\For{$k=0$ to $M-1$}
\State Given  $\cV^{k}$, choose $N_c$ virtual collision pairs $(i_\ell,i_{\ell_1})$ uniformly without replacement. The remaining $N-2N_c$ particles do not have a virtual (or real) collision and set $v^{k+1}_{i} = v^{k}_{i}$.
\For{$\ell= 1$ to $N_c$}
\State Compute $u_\ell = |v^k_{i_{\ell}}-v^k_{i_{\ell_1}}|^\beta$.
\State Draw a random number $\xi_\ell$ from the uniform distribution $\mathcal{U}([0,1])$.
\If{ $\xi_\ell \leq u_\ell/\Sigma_v $}
\State Sample scattering angle $\theta \sim C_\kappa (\theta)\sin(\theta)$ and azimuthal angle $\varphi\sim\mathcal{U}([0,2\pi])$.
\State Perform real collision between $ v^k_{i_{\ell }}$ and $v^k_{i_{\ell_1 }} $ following~\eqref{BoltzmannSolution} and obtain $({v^k_{i_{\ell }}}', {v^k_{i_{\ell_1}}}')$.
\State Set $(v^{k+1}_{i_{\ell }} , v^{k+1}_{i_{\ell_1}} ) = ({v^k_{i_{\ell }}}', {v^k_{i_{\ell_1}}}')$.
\Else 
\State The virtual collision is not a real collision. Set $( v^{k+1}_{i_{\ell }} , v^{k+1}_{i_{\ell_1}} ) = ({v^k_{i_{\ell }}}, {v^k_{i_{\ell_1}}}$).
\EndIf
\EndFor
\EndFor
\end{algorithmic}
\end{algorithm}

The DSMC method presented in~\Cref{alg:VHS_DSMC} applies to any general collision kernel $q(v-v_1, \sigma)$ with an upper bound $\Sigma$. However, many common collision models for the Boltzmann equation are of the form~\eqref{eq:VSS/VHS}, for which the collision kernel contains two separate parts: the relative velocity part $|v-v_1|^\beta$, and the scattering angle-dependent part $C_\kappa(\theta)$. Given this separation of variables, the DSMC method presented in~\Cref{alg:VHS_DSMC} can be modified to improve sampling efficiency. 

In contrast to the upper bound $\Sigma$ for the entire collision kernel as defined in~\eqref{eq:UpperBound}, we define $\Sigma_v$ as the upper bound for the relative velocity part where
\begin{equation}\label{eq:UpperBound_V}
    |v-v_1|^\beta < \Sigma_v.
\end{equation}
Numerically, the upper bound can be taken as $\Sigma_v = \max_i |v_i - v_{i_1}|^\beta$ where $i$ is the index of a velocity particle and $i_1$ is its collision partner index~\cite{pareschi2001introduction}. Besides, we define the weighted surface area
\begin{equation}\label{eq:beta Sphere area}
    A_\kappa = \int_{0}^{2 \pi} \int_{0}^{\pi} C_\kappa(\theta) \sin (\theta) d\theta d\varphi,
\end{equation}
which is different from the unweighted surface area $A$ defined in~\eqref{eq:mu}. Correspondingly, we have a different collision rate $\mu_\kappa = A_\kappa \Sigma_v \rho$ where the density $\rho = \int f dv$. We summarize the modified DSMC algorithm in~\Cref{alg:VHS_DSMC_2}.

Compared with~\Cref{alg:VHS_DSMC}, there is a major change in~\Cref{alg:VHS_DSMC_2}. The collision scattering angle $\theta$ and the azimuthal angle $\varphi$ are sampled independently and directly instead of using the rejection sampling. Since $\theta \in \mathbb R$, it is quite efficient to use inverse transform sampling to sample $C_\kappa(\theta)  \sin (\theta)$.

Note that~\Cref{alg:VHS_DSMC_2} can be used to sample angle-independent kernels, i.e., $q(v-v_1,\sigma) = |v-v_1|^\beta$. However, since the scattering and azimuthal angles are sampled uniformly first, there is a numerical dependence of the resulting $\sigma$ on $u = v-v_1$ based on~\cref{eq:sigma}, despite that there is no $\sigma$-dependence in the collision kernel itself. Technically, the adjoint Jacobian matrix $D$ should follow the formulation~\eqref{eq:new D} in the corresponding adjoint DSMC~\Cref{alg:adjoint_VHS} instead of~\eqref{eq:D angle-independent}. Nevertheless, we observe that the additional term $\widetilde B$ in~\eqref{eq:new D} can be dropped due to an averaging effect, reducing it to~\eqref{eq:D angle-independent}. The resulting gradient is not affected by $\widetilde B$. We will further demonstrate it numerically in~\Cref{subsec:test_angle_independent}.

On the other hand, \Cref{alg:VHS_DSMC} does not require sampling the scattering angle and the azimuthal angle first for a general angle-dependent kernel $q(v-v_1,\sigma)$. It samples $(v,v_1,\sigma)$ independently  before going through the rejection sampling step. Although $\sigma$ is originally sampled independent of $(v,v_1)$, its acceptance probability, $q(v-v_1,\sigma)$, enforces the dependence of the final accepted $\sigma$ on $(v,v_1)$. As a result, the adjoint Jacobian matrix $D$ in the corresponding adjoint DSMC~\Cref{alg:adjoint_VHS} should still have the $\widetilde B$ term defined in~\eqref{eq:tilde B def} when the forward DSMC method follows~\Cref{alg:VHS_DSMC}.


\section{Numerical Results}\label{sec:numerics}

In this section, we numerically test the adjoint DSMC method described in~\Cref{alg:adjoint_VHS} by computing the gradient using formula \eqref{eq:adjoint J^{AD}} for the objective function discussed in~\Cref{eq:OTD_obj,eq:obj1},
$$
J(m)= \int_{\bbR^3} \phi (v) f(v,T) dv  \approx \bar{\phi}^M = \frac{\rho}{N}\sum_{i=1}^N \phi(v_i^M),
$$
evaluated at the final time $t=T$, with respect to the parameter $m$ in the initial conditions $f_0(v;m)$ for the general collision kernel in the form of~\eqref{eq:VSS/VHS}, 
\begin{equation} \label{eq:C_beta}
q(v-v_1,\sigma) = C_\kappa(\theta) |v-v_1|^\beta =  \frac{1+\kappa}{2^{\kappa+2}\pi \epsilon}  (1+\cos{\theta})^\kappa |v-v_1|^\beta,  \, \,  \, \cos{\theta} = \sigma \cdot \frac{v-v_1}{|v-v_1|} .
\end{equation}
We choose this particular form of $C_\kappa(\theta)$ in~\eqref{eq:C_beta} such that $A_\kappa=1/\epsilon$ as defined in \eqref{eq:beta Sphere area}, where $\epsilon$ controls the amount of collisions per unit of time. 
We first consider collision kernels with only velocity dependence ($\kappa=0$ and $\beta=\{0,1,2\}$) and next we consider collision kernels with both velocity and angle dependence ($\kappa=\{1,2,5\}$ and $\beta=1$).
We also assume that  $\rho(t)=\int_{\bbR^3} f(v,t) dv=1$ (which is conserved throughout the evolution by the forward DSMC~\Cref{alg:VHS_DSMC,alg:VHS_DSMC_2}), and thus $\mu_\kappa = A_\kappa \Sigma_v \rho = \Sigma_v/\epsilon$.

Since we do not have an exact solution to the Boltzmann equation~\eqref{eq:homoBoltz} under this general collision kernel, we compare the gradient  $\nabla_m J^{AD}$ computed by the adjoint DSMC method (via \Cref{alg:adjoint_VHS} and \cref{eq:adjoint J^{AD}}) with the one computed by the central finite difference method,
\begin{equation}~\label{eq:FD_gradient}
\nabla_m J^{FD}(m) \approx \frac{J(m+\Delta m)-J(m-\Delta m)}{2 \Delta m},
\end{equation}
where both $J(m+\Delta m)$ and $J(m-\Delta m)$ are computed using forward DSMC simulations with the same random seed. The random seed is used to initialize a pseudorandom number generator for sampling steps in the forward DSMC algorithms (see~\Cref{alg:VHS_DSMC,alg:VHS_DSMC_2}). We fix the random seed in~\eqref{eq:FD_gradient} to reduce the random error in approximating $\nabla_m J^{FD}$. We will also use the average value of $\nabla_m J^{FD}$ from multiple runs (under different random seeds) to further reduce the variance in $\nabla_m J^{FD}$.

We are interested in the error between the two gradients defined as $\big|\nabla_m J^{AD} - \nabla_m J^{FD}\big|$, and studying its convergence as the number of particles $N$ increases in the empirical distribution~\eqref{eq:empirical} as a discretization for the distribution function $f$ in the Boltzmann equation~\eqref{eq:homoBoltz}. 
This total error has several contributions: the random error from particle discretization in both $\nabla_m J^{AD}$ and $\nabla_m J^{FD}$, and the finite difference error in $\nabla_m J^{FD}$. Thus, when we study the convergence of the total error with $N$ increasing, we expect the total error first to decrease and then to plateau at a fixed constant determined by the finite difference error after the random error is sufficiently small. 

To further reduce the random error in both $\nabla_m J^{AD}$ and $\nabla_m J^{FD}$, we perform $M_s=100$ computations for each of them using random initial conditions, and compute the average values, denoted as $\overline{\nabla_m J}^{AD}$ and $\overline{\nabla_m J}^{FD}$,  before computing the error $e$.  That is,
\begin{equation}~\label{eq:error}
e=\big|\overline{ \nabla_m J}^{AD} - \overline{ \nabla_m J}^{FD}\big|.
\end{equation}
The random errors in the averaged gradients are estimated by first computing the standard deviation in $\nabla_m J^{AD}$ and $\nabla_m J^{FD}$ respectively using the $M_s$ i.i.d.~runs,  followed by a rescaling using the factor $1/\sqrt{M_s}$.  

Ultimately, we want to show that $| \overline{\nabla_m J}^{AD} - \nabla_m J|$ is small where $\nabla_m J$ is the true gradient.  To make $\overline{\nabla_m J}^{FD}$ a better approximation of $\nabla_m J$,  one needs to use a tiny perturbation $\Delta m$ to reduce the central difference error of $\mathcal{O}(|\Delta m|^2)$, and relatively large $N$ and $M_s$ to reduce the random error of $\mathcal{O}((\sqrt{N M_s} |\Delta m|)^{-1} )$.  In our previous work~\cite{caflisch2021adjoint}, we balanced the particle discretization (random) error and the finite difference error in $\overline{\nabla_m J}^{FD}$ and found that using $\Delta m=0.1$ was nearly optimal for $N=10^6$ Maxwellian particles, $M_s=100$, and initial conditions similar to \eqref{eq:IC}. Hence, we will also fix $\Delta m=0.1$ for all simulations here as a heuristic estimation in the following numerical tests.

For the function $\phi(v)$ in \eqref{eq:OTD_obj}, we use $v_l^2$, $l\in\{x,y,z\},$ so the objective functions are 
$$T_l=\frac{1}{N}\sum_{i=1}^N{(v_{F,i}^l)^2 } \approx \int_{\bbR^3} v_l^2 f(v,T) dv,\quad l\in\{x,y,z\},$$ the second-order velocity moments of the distribution function in the $l$-direction at the time $t=T$. For the parameter $m$, we use temperature values in the initial distribution function $m=(T_x^0, T_y^0, T_z^0)$. We further refer to these gradients as $\frac{\delta J}{\delta m}=\frac{\partial T_l}{\partial T_p^0}$, $l,p\in\{x,y,z\}$. 
Here, the dimension of the parameter $m$ is only $3$. A real advantage of this adjoint-state method is, of course, when the parameter $m$ is extremely high-dimensional since the adjoint-state method allows one to compute all the components of the gradient $\nabla_m J^{AD}$ by doing only \textit{one} forward DSMC simulation and \textit{one} backward adjoint DSMC simulation.

In all the numerical tests below, we use the same anisotropic Gaussian as the initial condition,
\begin{equation}~\label{eq:IC}
f_0(v) = \frac{1}{(2\pi)^{3/2}\sqrt{T_x^0 T_y^0 T_z^0}} \exp\left(-\frac{v_x^2}{2T_x^0}-\frac{v_y^2}{2T_y^0}-\frac{v_z^2}{2T_z^0}\right),
\end{equation}
where $T_x^0=1,T_y^0=1,T_z^0=0.5$, so that $m=(1,1,0.5)$. In this case, the solution to the Boltzmann equation~\eqref{eq:homoBoltz} will relax to an isotropic Gaussian with the temperature $T_M=(T_x^0+T_y^0+T_z^0)/3=0.8333(3)$ as time increases. In the following tests, we use~\Cref{alg:VHS_DSMC_2} for the forward DSMC simulations with a time step $\Delta t=0.1$ and the final time $t=T=M \Delta t=2$ at which only partial relaxation to the isotropic Gaussian is attained. Since the initial condition and the solution are anisotropic Gaussians with $T_x^0,T_y^0,T_z^0 \leq 1$, then the upper bound for the relative velocity, $\max_i |v_i - v_{i_1}|$, can be taken to be $10$, and consequently we can set $\Sigma_v =10^\beta$ in \eqref{eq:UpperBound_V}. Therefore, the fraction of particles that participate in virtual collisions at every time step of the forward DSMC simulations is $ N_c/N = \Delta t \mu_\kappa = \Delta t \Sigma_v / \epsilon = \Delta t 10^\beta / \epsilon$,  based on~\Cref{alg:VHS_DSMC_2}.  As discussed in~\Cref{sec:algs13_discuss}, for collision kernels in the form of~\eqref{eq:VSS/VHS}, it is more computationally efficient to use~\Cref{alg:VHS_DSMC_2} as the forward DSMC algorithm than~\Cref{alg:VHS_DSMC}.  In the numerical tests below, we will also only use~\Cref{alg:VHS_DSMC_2}, but will comment on the gradient calculation based on~\Cref{alg:VHS_DSMC} in~\Cref{sec:algs13_comments}.

To compute the gradient $\nabla_m J^{AD}$ using Algorithm \ref{alg:adjoint_VHS} and formula \eqref{eq:adjoint J^{AD}}, we need $\nabla_m  v_{i}^0 $. Since our initial distribution \eqref{eq:IC} is an isotropic Gaussian, we can sample $v_{i}^0 $ by sampling $3N$ values from the standard normal distribution $\cN(0,1)$ and then rescaling the values with appropriate initial temperatures as 
$$v_i^0=(v_i^{x,0},v_i^{y,0},v_i^{z,0}) = (\sqrt{T_x^0}\cN_i^{x,0} ,\sqrt{T_y^0}\cN_i^{y,0},\sqrt{T_z^0}\cN_i^{z,0}),$$ 
where $\cN_i^{x,0},\cN_i^{y,0},\cN_i^{y,0}$ are samples of $\cN(0,1)$. For the parameter $m=T_p^0$, we can then compute 
$$
\nabla_m  v_i^{j,0}= \frac{ \cN_i^{j,0} }{2\sqrt{T_p^0}} \delta_{j,p}=\frac{ v_i^{j,0}}{2T_p^0} \delta_{j,p},\quad j,p\in\{x,y,z\}.$$
This way of obtaining $\nabla_m  v_i^{j,0}$ corresponds to the so-called pathwise gradient estimator~\cite[Section 5]{mohamed2019monte}. 

\begin{figure}
     \centering
    \subfloat[Error $e$ computed without~\eqref{eq:tilde B def}]{\includegraphics[width=.5\linewidth]{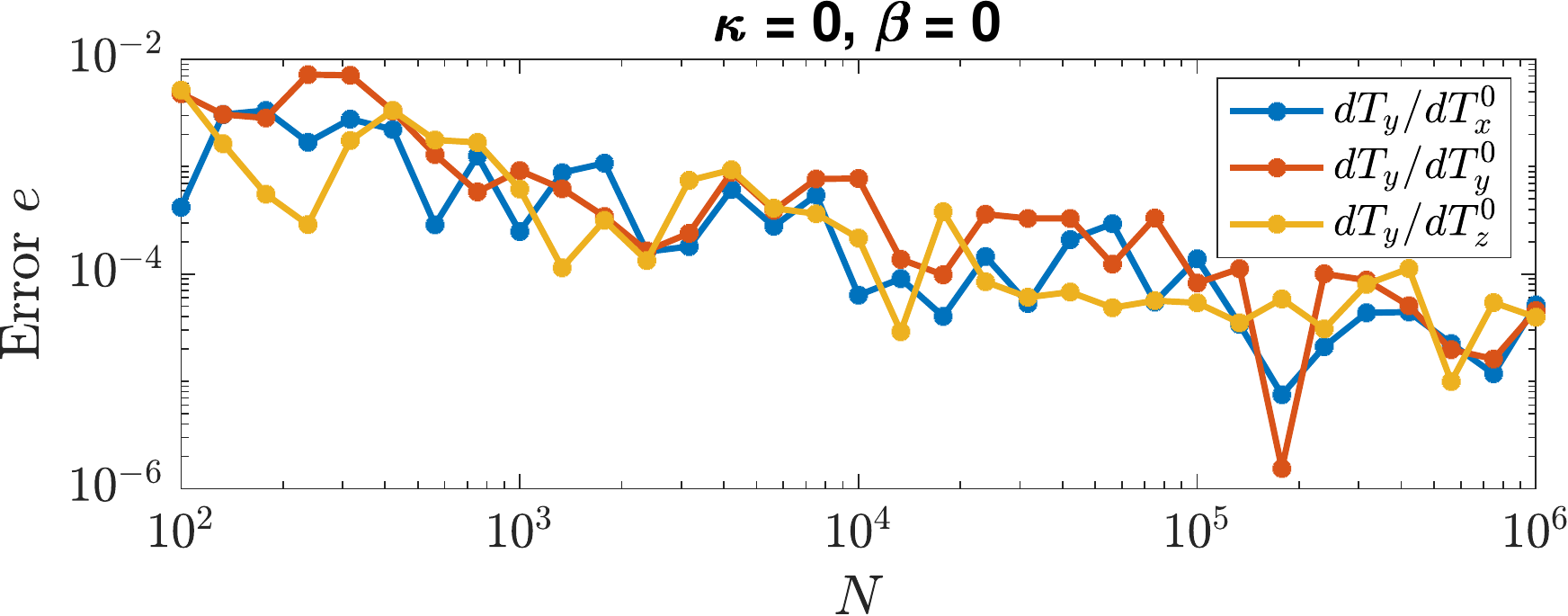}\label{fig:Max-wo}}
     \subfloat[Error $e$ computed with~\eqref{eq:tilde B def}]{\includegraphics[width=.5\linewidth]{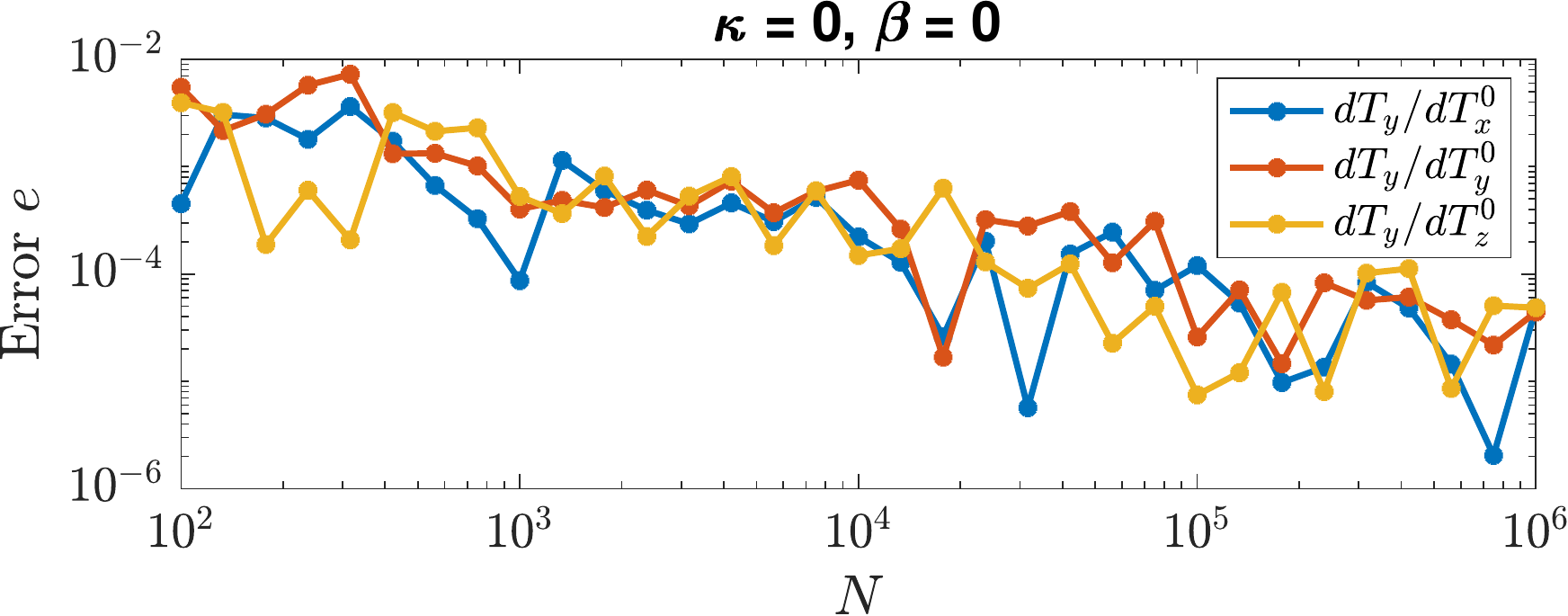}\label{fig:Max-wt}}\\
        \subfloat[Differences in $\overline{\nabla_m J}^{AD}$ with and without~\eqref{eq:tilde B def}]{\includegraphics[width=.5\linewidth]{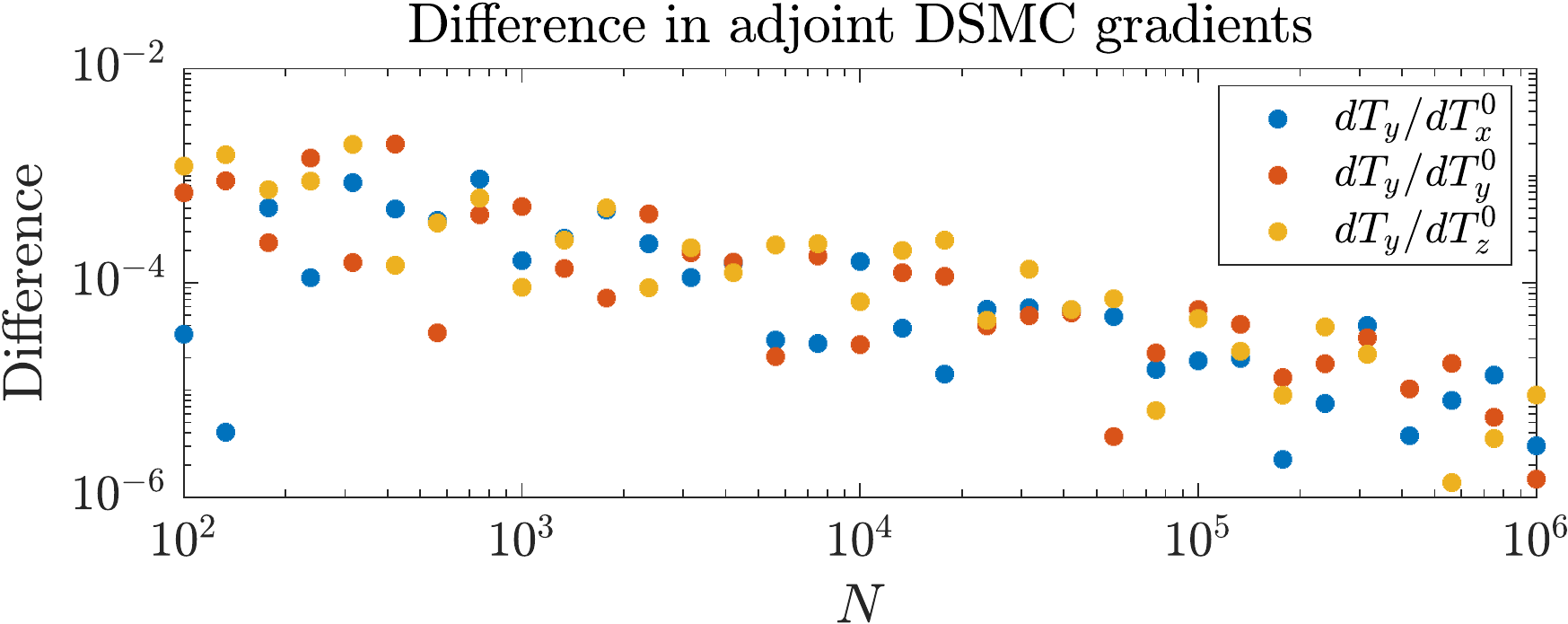}\label{fig:Max-diff}}
     \subfloat[The random errors in $\overline{\nabla_m J}^{AD}$]{\includegraphics[width=.5\linewidth]{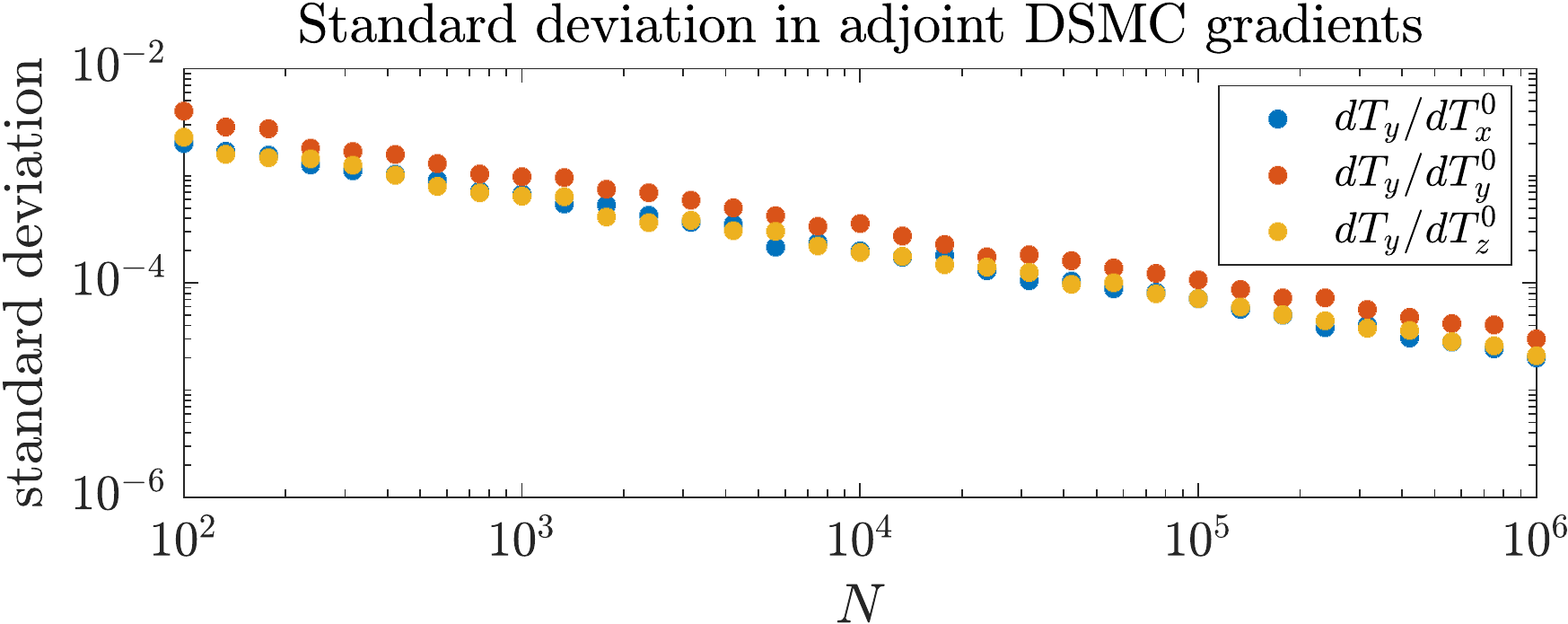}\label{fig:Max-std}}\\  
     \caption{(a)-(b):~Comparison of the gradient error~\eqref{eq:error} for the Maxwellian constant collision kernel ($\kappa = 0$, $\beta = 0$) with and without the term~\eqref{eq:tilde B def} in the adjoint matrix~\eqref{eq:new D} in the adjoint DSMC~\Cref{alg:adjoint_VHS}.  
(c):~The absolute value of the difference in $\overline{\nabla_m J}^{AD}$ between the two cases. (d):~The random error in $\overline{\nabla_m J}^{AD}$ measured by the estimated standard deviation in the averaged value $\overline{\nabla_m J}^{AD}$ from $M_s = 100$ i.i.d.~runs. Note that the term~\eqref{eq:tilde B def} is deterministic and thus does not affect the random error in the adjoint DSMC gradient.}
     \label{fig:Max_wtwo_compare}
\end{figure}

\subsection{Simulations with Angle-Independent Collision Kernels}\label{subsec:test_angle_independent}
First,  we focus on the angle-independent collision kernel by setting $\kappa=0$ and $\epsilon=10$. Thus,  the collision kernel $q(v-v_1,\sigma) = 1/(40 \pi) |v-v_1|^\beta$ following~\eqref{eq:C_beta}.  

When $\beta = 0$,  the collision kernel corresponds to Maxwell molecules discussed in~\cite{caflisch2021adjoint}.  We focus on the objective function $T_y$ as an example. After $M=20$ time steps of the forward DSMC simulation using different numbers of particles $N$,  we numerically illustrate the error in the gradient calculation both with and without the term~\eqref{eq:tilde B def} in the adjoint equation; see~\Cref{fig:Max-wo,fig:Max-wt}.  The differences in the averaged adjoint DSMC gradients $\overline{\nabla_m J}^{AD}$, as seen in~\Cref{fig:Max-diff}, are observed to be in the same order as the random errors in $\overline{\nabla_m J}^{AD}$, as shown in~\Cref{fig:Max-std}.  The comparison in~\Cref{fig:Max_wtwo_compare} illustrates that the additional term does not affect angle-independent kernels due to an averaging effect.   This phenomenon occurs because, for angle-independent collision kernels, the post-collision relative velocity $\sigma$ can be seen as uniformly distributed over the sphere, which weakens its dependence on the scattering angle $\theta$ and the pre-collision relative velocity $\alpha$, despite the relation $\cos \theta = \sigma \cdot \alpha$. In~\cite{caflisch2021adjoint},  we assumed that $\sigma$ does not depend on $\alpha$ for the Maxwellian gas and obtained the adjoint DSMC algorithm without the term~\eqref{eq:tilde B def}. Here, we use this example to illustrate that both adjoint matrices are valid for angle-independent kernels.

\begin{figure}
     \centering
     \subfloat[Gradient error $e$ in~\eqref{eq:error}]{\includegraphics[width=.46\linewidth]{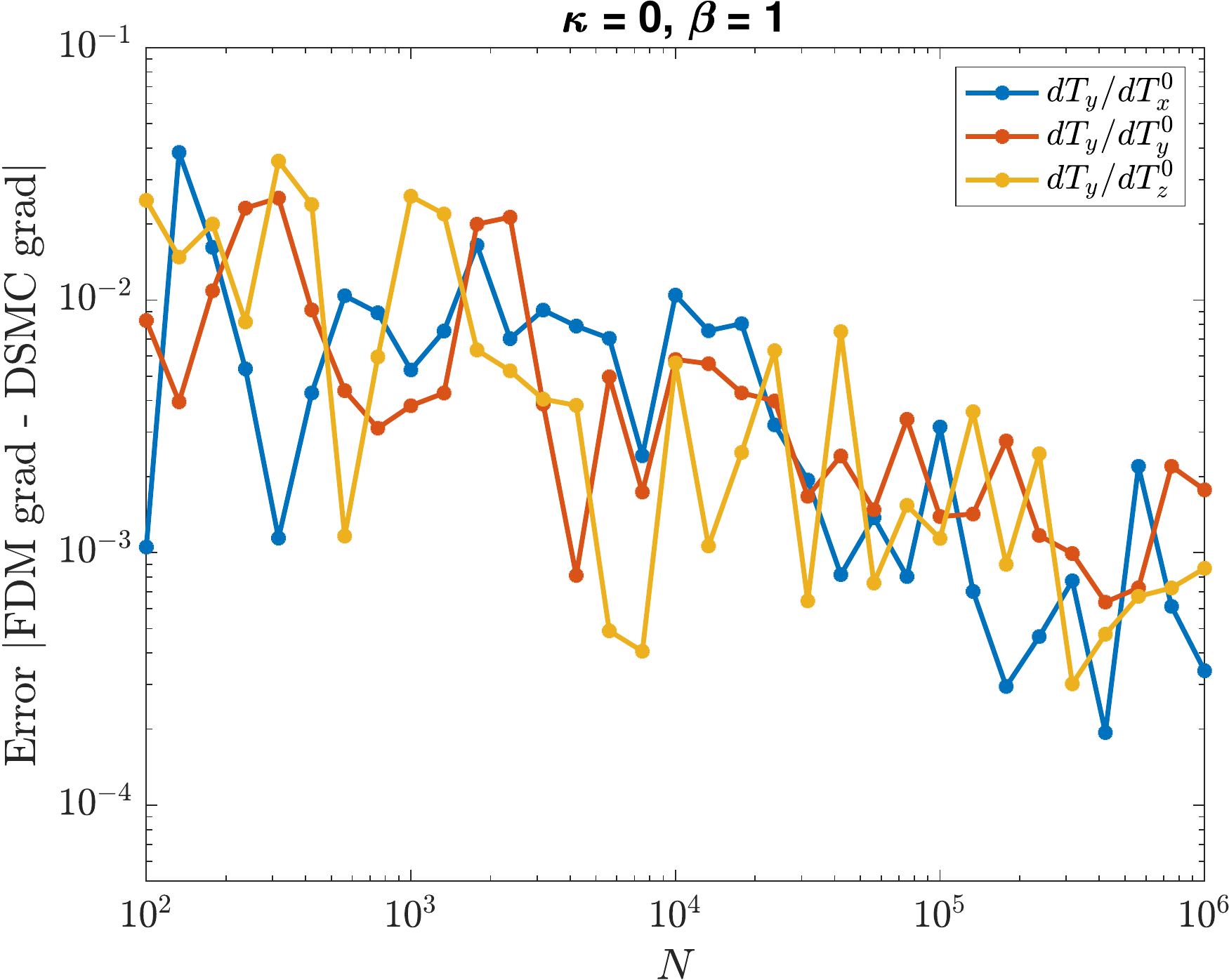}\label{fig:beta1_error}}
     \hspace{0.07\linewidth}
     \subfloat[Standard deviations of $\overline{\nabla_m J}^{FD}$ and $\overline{\nabla_m J}^{AD}$]{\includegraphics[width=.46\linewidth]{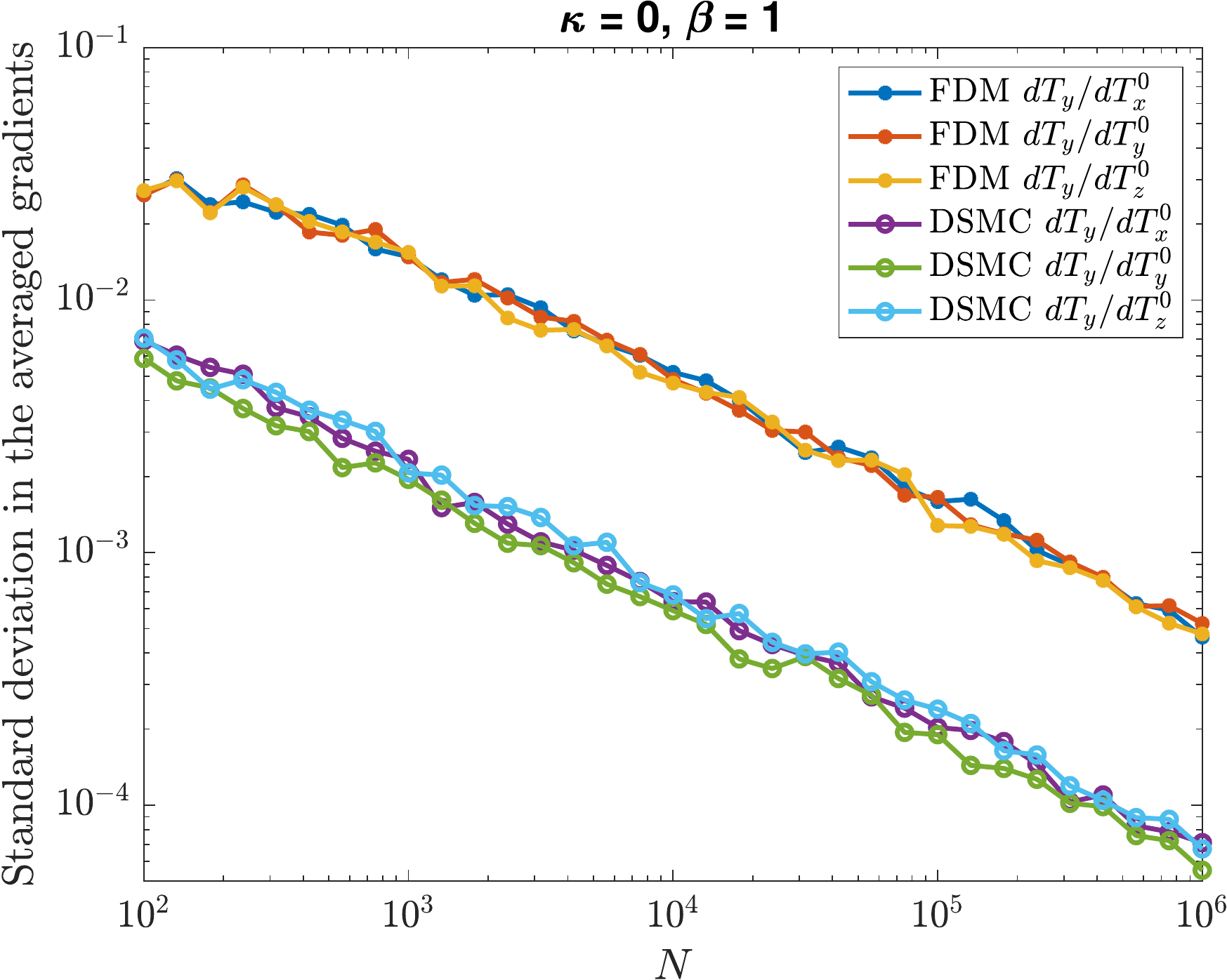}\label{fig:beta1_std}}
     \caption{The gradient error (a) and the standard deviation (b) for the case $\beta=1$ and $\kappa=0$ after $20$  time steps.  The standard deviations plotted in (b) are for the averaged gradient values from $M_s = 100$ i.i.d.~runs. }
\end{figure}

Next, we consider the case $\beta=1$.  The objective function is $T_y$, the final temperature in the $y$ direction. In~\Cref{fig:beta1_error}, we show that the gradient error $e$, defined in~\eqref{eq:error}, decays as the number of particles $N$ increases,  while~\Cref{fig:beta1_std} illustrates that the standard deviations of the averaged adjoint DSMC gradient $\overline{\nabla_m J}^{AD}$ and finite-difference gradient $\overline{\nabla_m J}^{FD}$ both decay as $\mathcal{O}(1/\sqrt{N})$.  As mentioned earlier, the error $e$ observed in~\Cref{fig:beta1_error} has three contributions: the random error contributions from both $\nabla_m J^{AD}$ and $\nabla_m J^{FD}$, and the finite difference error in $\nabla_m J^{FD}$.  We remark that, based on the standard deviations shown in~\Cref{fig:beta1_std},  the random error in ${\nabla_m J}^{FD}$ could be a leading contribution in $e$ and the random error in the adjoint DSMC gradient ${\nabla_m J}^{AD}$ is nearly $10$ times smaller.  



\begin{figure}
     \centering
     \subfloat[Error $e$ defined in~\eqref{eq:error}]{\includegraphics[width=.7\linewidth]{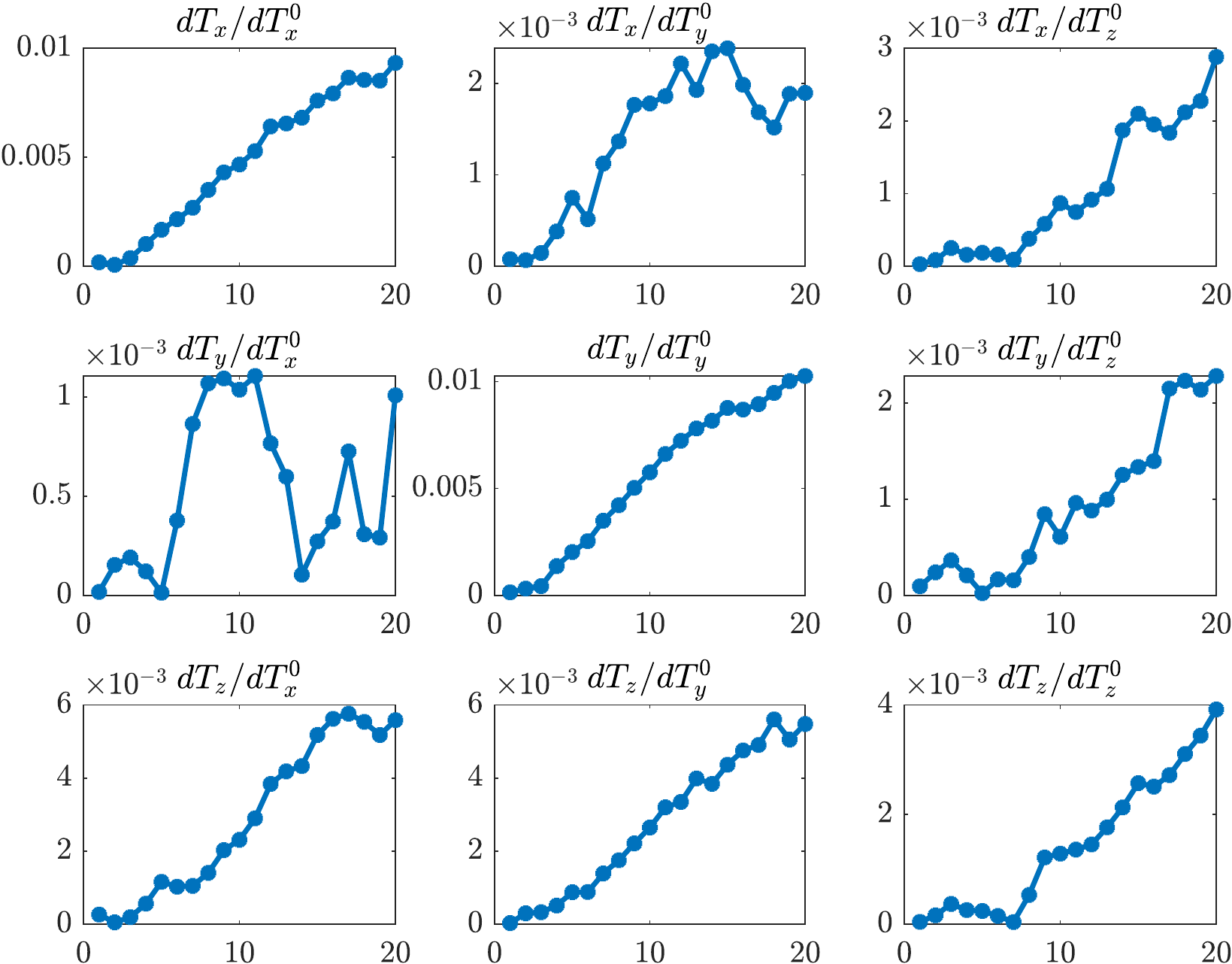}\label{fig:beta2 error}}\\
     \vspace{1cm}
     \subfloat[Standard deviations in the averaged gradients $\overline{\nabla_m J}^{FD}$ and $\overline{\nabla_m J}^{AD}$]{\includegraphics[width=.7\linewidth]{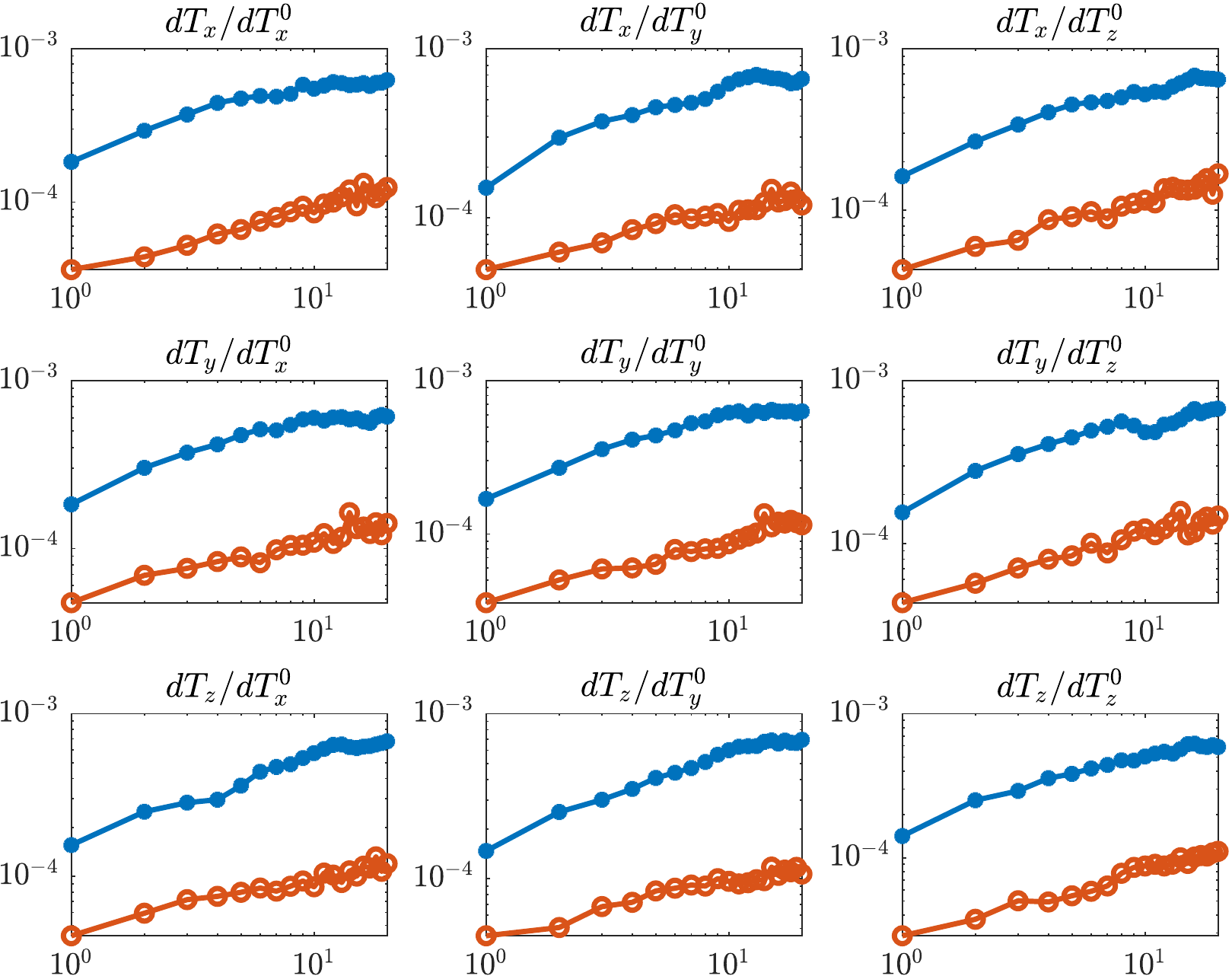}\label{fig:beta2 std}}
     \caption{The gradient error (a) and the standard deviation (b) of $\frac{\partial T_l}{ \partial T^0_p}$, $l,p\in\{x,y,z\}$, for the case $\beta=2$, $\kappa=0$ and $N = 10^6$,  with the number of time steps $M$ ranging from $1$ to $20$ (the $x$ axis).  In~\Cref{fig:beta2 std},  the blue log-log plots represent the standard deviations of $\overline{\nabla_m J}^{FD}$ while the red log-log plots are for $\overline{\nabla_m J}^{AD}$.  }
     \label{fig:beta2}
\end{figure}

In~\Cref{fig:beta2}, we consider the case $\beta=2$, $\kappa=0$,  and evaluate the  gradient $\frac{\partial T_l}{ \partial T^0_p}$,  where $l,p\in\{x,y,z\}$.  We fix the number of particles $N = 10^6$, and investigate how the gradient error $e$ defined in~\eqref{eq:error} and the standard deviations of $\overline{\nabla_m J}^{AD}$ and $\overline{\nabla_m J}^{FD}$ change with respect to the total number of time steps $M$  (the $x$ axis in the plots).  We consider $M = 1,\ldots, 20$.   In~\Cref{fig:beta2 error},  the gradient errors are mostly linearly increasing up to perturbations incurred by the random errors (illustrated in~\Cref{fig:beta2 std}). Based on the forward DSMC~\Cref{alg:VHS_DSMC,alg:VHS_DSMC_2}, the variations of both gradients,  $\overline{\nabla_m J}^{AD}$ and $\overline{\nabla_m J}^{FD}$,  are expected to be $\mathcal{O}(M)$ since the number of sampling steps in the forward DSMC algorithms grows linearly in time.  Therefore,  we observe from the log-log plots in~\Cref{fig:beta2 std} that the standard deviations for both gradients grow as $\mathcal{O}(\sqrt{M})$, while the standard deviation of the finite difference gradient  $\overline{\nabla_m J}^{FD}$ is much bigger than the one for the adjoint DSMC gradient $\overline{\nabla_m J}^{AD}$.

%
%
%
%

\subsection{Simulations with Angle-Dependent Collision Kernels}

\begin{figure}
     \centering
     \subfloat[Gradient error $e$ for $\kappa = 1$,   $\beta = 1$]{\includegraphics[width=.45\linewidth]{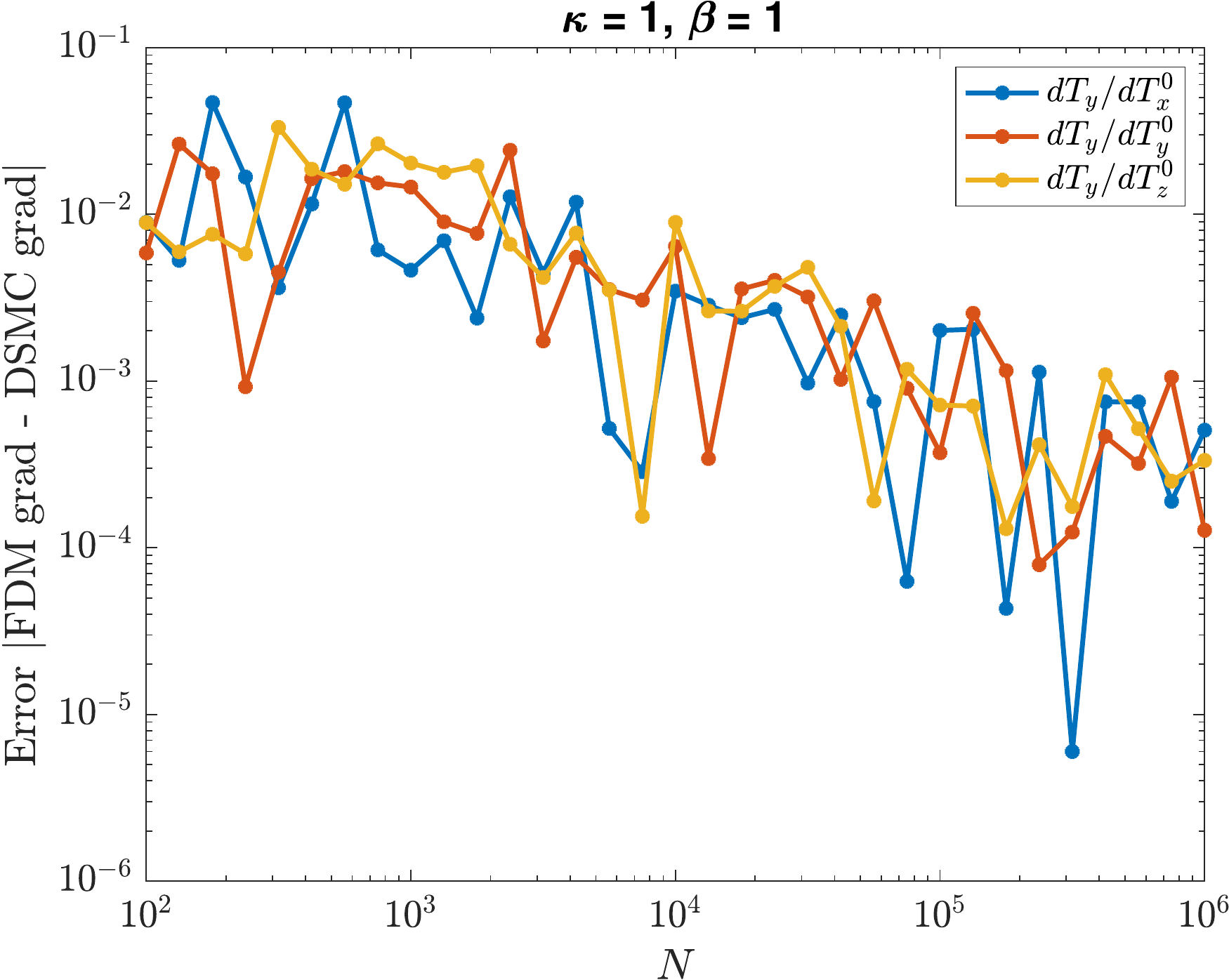}\label{fig:Kappa1 error}}
     \hspace{0.08\textwidth}
     \subfloat[Gradient error $e$ for $\kappa = 2$,  $\beta = 1$]{\includegraphics[width=.45\linewidth]{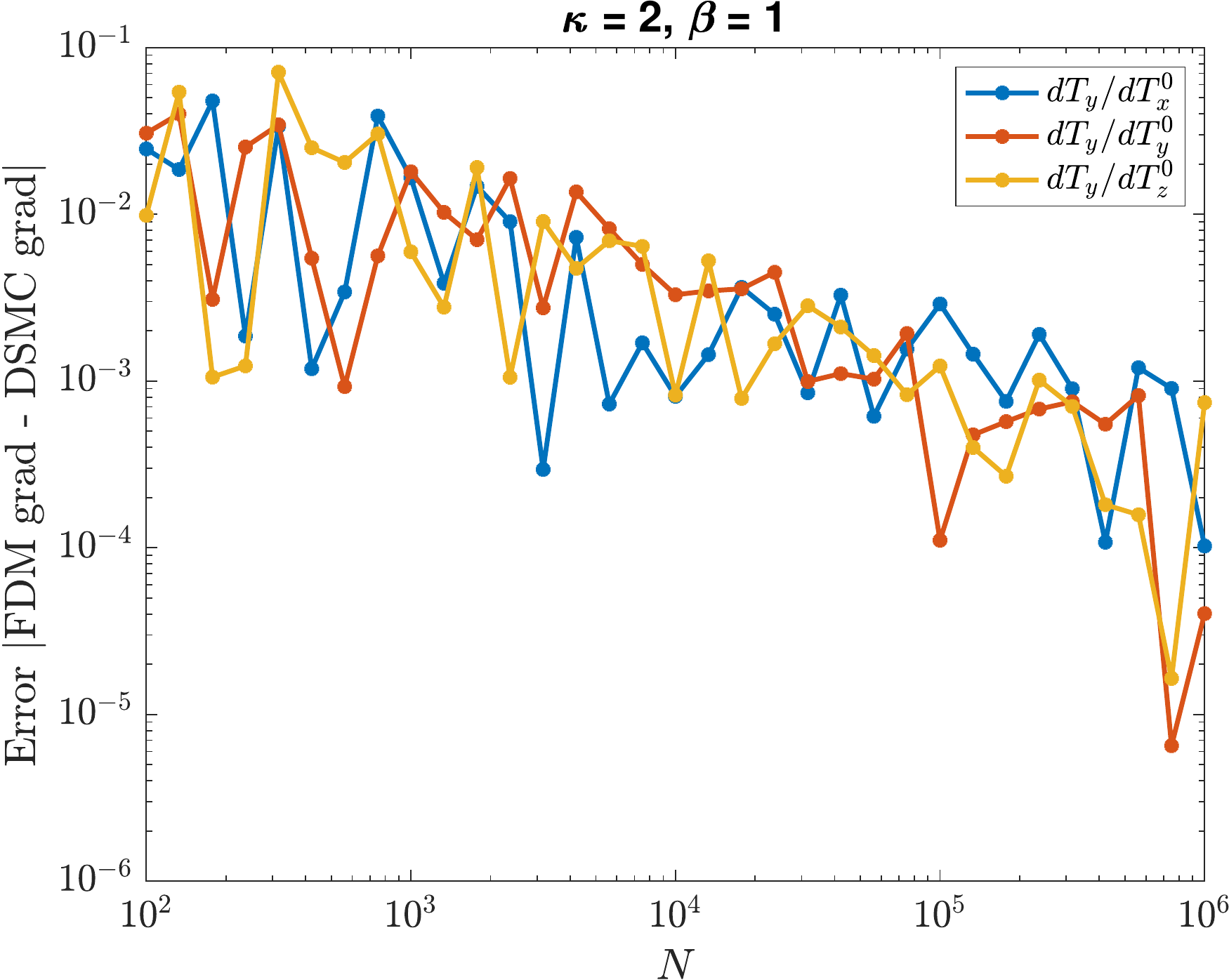}\label{fig:Kappa2 error}}
     \caption{The gradient error~\eqref{eq:error} (the $y$ axis) for the collision kernel~\eqref{eq:C_beta} with $\beta = 1$, $\kappa = 1$ (left) and $\kappa = 2$ (right),  after $M=20$ time steps.  The number of particles $N$ (the $x$ axis) ranges from $10^2$ to $10^6$. }
     \label{fig:kappa12}
\end{figure}

In this subsection, we consider collision kernels that are both velocity and angle-dependent. We set $\beta=1$,  $\epsilon=10$ and consider various $\kappa$ values.   The collision kernel takes the form 
$$q(v-v_1,\sigma) = \frac{ (1+\kappa) }{10 \pi \, 2^{\kappa+2} }(1+\cos{\theta})^\kappa |v-v_1|,\quad \cos \theta  = \sigma \cdot \frac{v-v_1}{|v-v_1|}.$$
We remark that when the collision kernel is angle-dependent, the adjoint matrix $D_i^k$ in~\Cref{alg:adjoint_VHS} should follow eq.~\eqref{eq:new D} instead of~\cref{eq:D angle-independent}.  Note that the difference between the two adjoint matrices is that eq.~\eqref{eq:new D} has an additional term~\eqref{eq:tilde B def}.

\begin{figure}
     \centering
     \subfloat[Error $e$ computed without~\eqref{eq:tilde B def}]{\includegraphics[width=.7\linewidth]{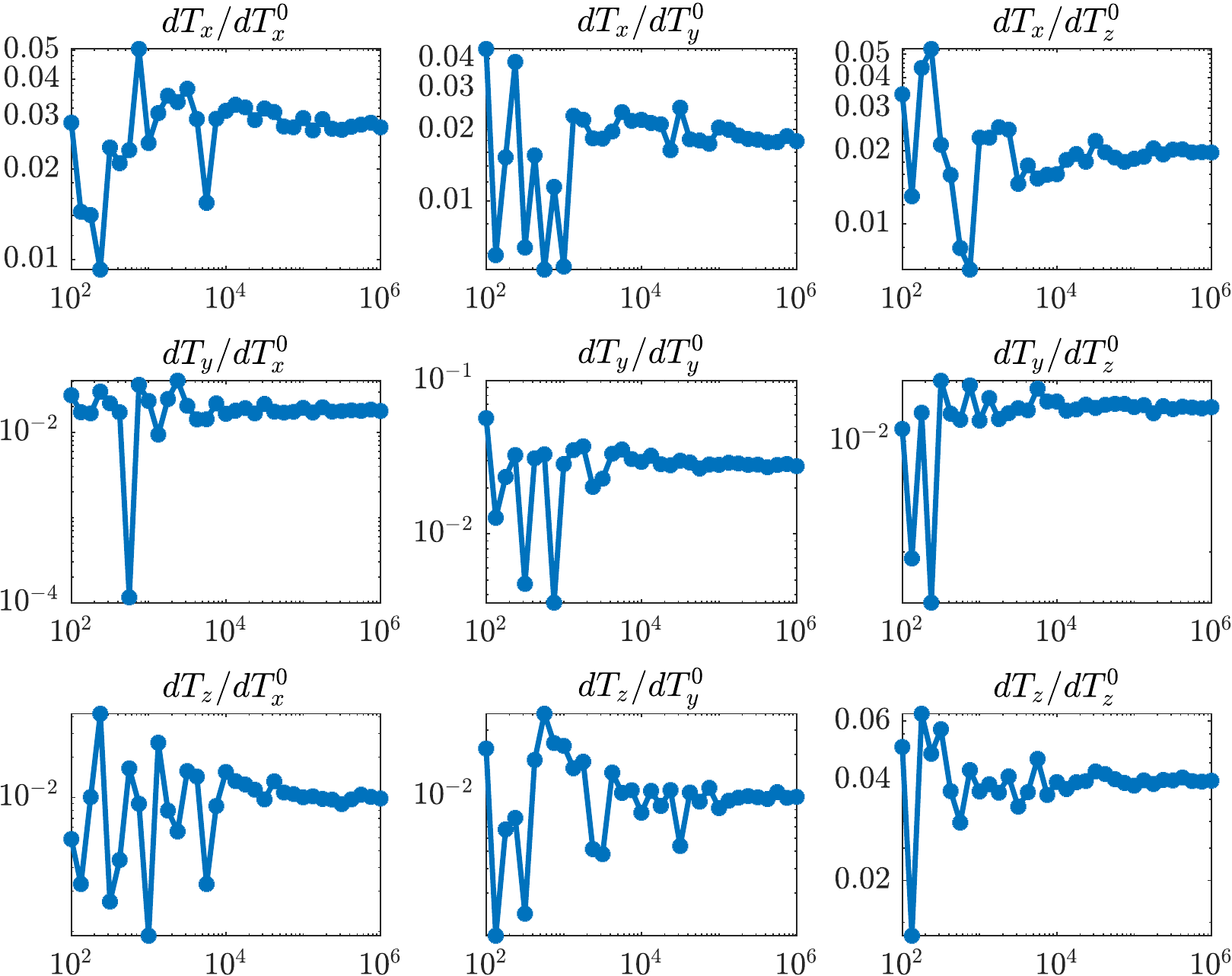}\label{fig:Kappa5without}}\\
     \vspace{1cm}
          \subfloat[Error $e$ computed with~\eqref{eq:tilde B def}]{\includegraphics[width=.7\linewidth]{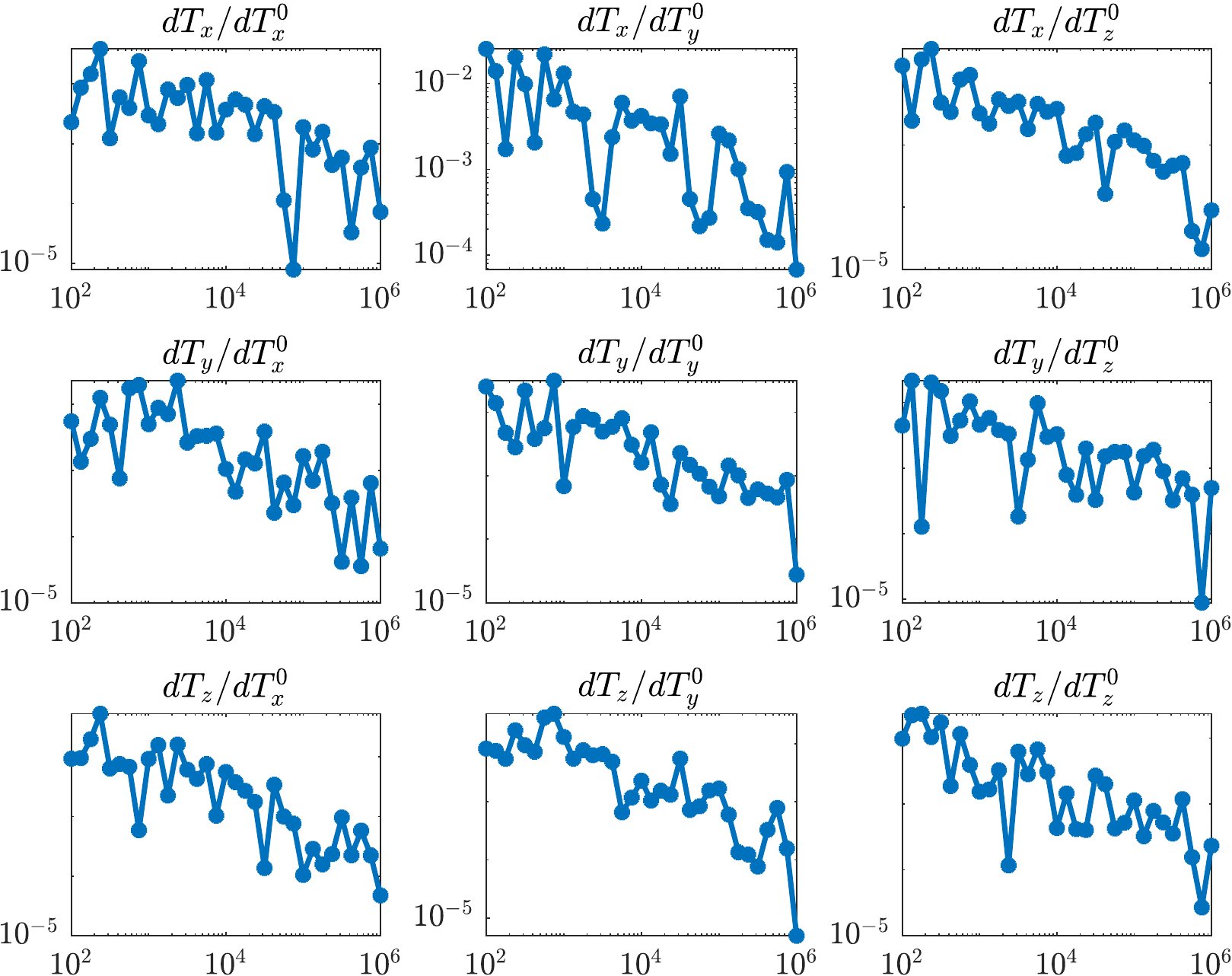}\label{fig:Kappa5with}}
     \caption{Comparison of the gradient error~\eqref{eq:error} (the $y$ axis) for the collision kernel with $\kappa = 5$, $\beta = 1$ after $M=20$ time steps,  with and without the additional term~\eqref{eq:tilde B def} in the adjoint matrix~\eqref{eq:new D} in the adjoint DSMC~\Cref{alg:adjoint_VHS} for gradient calculation.  The number of particles $N$ (the $x$ axis) ranges from $10^2$ to $10^6$.  The gradient errors are large without~\eqref{eq:tilde B def} in the adjoint matrix,  which leads to a wrong adjoint equation.}
     \label{fig:kappa5}
\end{figure}

We  first consider cases $\kappa=1$ and $\kappa = 2$, and set the objective function as $T_y$ with the parameter $m= T^0_p$,  where $p \in \{x,y,z\}$.  We plot the gradient errors in~\Cref{fig:kappa12},  both of which decay as the number of particles $N$ increases. In~\Cref{fig:kappa5}, we focus on the case $\kappa = 5$ and evaluate the gradient $\frac{\partial T_l}{ \partial T^0_p}$ where $l,p\in\{x,y,z\}$, after $M=20$ time steps.   Similar to~\Cref{fig:Max_wtwo_compare},  we compare the gradient errors when the adjoint DSMC gradients $\nabla_m J^{AD}$ are computed with the adjoint matrix~\eqref{eq:D angle-independent} and~\eqref{eq:new D}, respectively.  In~\Cref{fig:Kappa5without}, all $9$ gradient errors plateaued at a relatively large constant after $N\geq 10^4$.  On the other hand,  in~\Cref{fig:Kappa5with},  all $9$ gradient errors asymptotically decay as the number of particles $N$ increases. All the gradient errors are about $1\times 10^{-5}$ when $N = 10^6$.  The finite difference gradients $\nabla_m J^{FD}$ are the same in both figures,  so the drastic differences in $e$ come from the adjoint DSMC gradients $\nabla_m J^{AD}$. The comparison in~\Cref{fig:kappa5} indicates that the adjoint matrix $D$ defined in~\eqref{eq:D angle-independent} is incorrect for angle-dependent kernels,  which leads to wrong adjoint DSMC gradients in~\Cref{fig:Kappa5without}. The difference between~\Cref{fig:Max_wtwo_compare} and~\Cref{fig:kappa5} also shows that the additional term~\eqref{eq:tilde B def} is crucial for angle-dependent kernels, but plays little role for angle-independent collision kernels.

\subsection{Comments on the General Kernel Case}\label{sec:algs13_comments}
All the tests above were performed for separable collision kernels of type \eqref{eq:VSS/VHS} and therefore were based on~\Cref{alg:VHS_DSMC_2} for the forward DSMC simulations and~\Cref{alg:adjoint_VHS} for the adjoint DSMC simulations. For general collision kernels $q(v-v_1,\sigma) =  \tilde q(|v-v_1|,\theta)$ that satisfy \eqref{eq:UpperBound}, we can use~\Cref{alg:VHS_DSMC} for the forward DSMC simulations together with~\Cref{alg:adjoint_VHS} for the adjoint DSMC simulations. To this end, we have also numerically verified an approach of computing the forward DSMC via~\Cref{alg:VHS_DSMC} (and corresponding adjoint DSMC via~\Cref{alg:adjoint_VHS}) that does not split the kernel~\eqref{eq:VSS/VHS} into velocity-dependent and angle-dependent parts, but rather performs acceptance-rejection samplings over the full collision kernel with the scattering angle $\sigma$ uniformly sampled over the unit sphere (as described in in~\Cref{alg:VHS_DSMC}).
In this case, one must modify the $\eta_i^k$ term in~\cref{eq:back_eta} used in~\Cref{alg:adjoint_VHS}. Previously, it was based on rejection sampling using $q \sim |v-v_1|^\beta$, whereas in this case it is based on $q \sim |v-v_1|^\beta (1+\cos \theta)^\kappa$. 

Based on our numerical results, we conclude that we still need the extra term \eqref{eq:tilde B def} in \eqref{eq:new D} when running~\Cref{alg:adjoint_VHS} for an angle-dependent collision kernel. Even though the collision parameters $\sigma$ for collisions are sampled uniformly in the forward DSMC~\Cref{alg:VHS_DSMC} and do not \textit{explicitly} depend on the collision velocity pair $(v,v_1)$, after the acceptance-rejection step, the unit vector $\sigma$ becomes \textit{implicitly} dependent on $(v,v_1)$. The only case where the term~\eqref{eq:tilde B def} is not needed in the adjoint matrix $D$ is when the collision kernel is angle-independent; see~\Cref{fig:Max_wtwo_compare} for an illustration. 

It is worth noting that using~\Cref{alg:VHS_DSMC} for kernel in the form of~\eqref{eq:VSS/VHS} is less efficient and takes more computational time compared to using~\Cref{alg:VHS_DSMC_2}. This is because it involves sampling more virtual collisions, and taking additional time to sample collision parameters $\sigma$ for all virtual collision pairs instead of sampling collision parameters only for the actual collisions (as done in~\Cref{alg:VHS_DSMC_2}). Nevertheless, this approach is more general and can be used for general collision kernels \eqref{alg:VHS_DSMC_2} as long as the condition \eqref{eq:UpperBound} is satisfied (numerically).


\section{Conclusions}\label{sec:conclusion}

As discussed in~\Cref{sec:Intro} and~\Cref{section:Influence}, the method developed in this work is mainly based on the DTO approach but also involves an OTD step. The reason for the OTD step is that the rejection sampling involves a decision of whether a virtual collision is real, and the decision depends on the unknown parameter for which we want to compute the gradient. Differentiation after this choice would be applied to a discontinuous function (to indicate whether a collision is real or virtual) of the particle velocity, leading to a singularity. The OTD step in the new adjoint DSMC method, i.e., first differentiating the expectation over the decision and then sampling the resulting gradient, enables us to circumvent the difficulty of directly differentiating a discontinuous decision function from the rejection sampling. 

Another main contribution of this work is to consider collision kernels that are also scattering angle-dependent. As a result, the post-collision relative velocity depends on the pre-collision relative velocity. Although in the angle-independent cases, such as a constant Maxwellian collision kernel, the post-collision relative velocity always depends on the scattering angle as well as the pre-collision velocities based on their definitions,  we can ignore this dependence due to the averaging effect since the post-collision relative velocity is uniformly distributed over the sphere. This is not the case for angle-dependent kernels. In our new derivations for the adjoint DSMC algorithm, the resulting adjoint equation has an additional term that reflects this dependence. We remark that this additional term can be efficiently computed without extra memory requirement. 

This paper extends the adjoint DSMC method to a much more general class of collision kernels for the Boltzmann equation than the original proposal~\cite{caflisch2021adjoint}.  In future works, we plan to extend the adjoint DSMC method, for example, to Coulomb collisions and apply the method to large-scale optimization problems constrained by Boltzmann equations.

\section*{Acknowledgements}
This material is based upon work supported by the National Science Foundation under Award Number DMS-1913129 and the U.S.~Department of Energy under Award Number DE-FG02-86ER53223. Y.~Yang acknowledges support from Dr.~Max R\"ossler, the Walter Haefner Foundation and the ETH Z\"urich Foundation.


\bibliographystyle{elsarticle-harv}
\bibliography{main}

\begin{thebibliography}{23}
\expandafter\ifx\csname natexlab\endcsname\relax\def\natexlab#1{#1}\fi
\providecommand{\url}[1]{\texttt{#1}}
\providecommand{\href}[2]{#2}
\providecommand{\path}[1]{#1}
\providecommand{\DOIprefix}{doi:}
\providecommand{\ArXivprefix}{arXiv:}
\providecommand{\URLprefix}{URL: }
\providecommand{\Pubmedprefix}{pmid:}
\providecommand{\doi}[1]{\href{http://dx.doi.org/#1}{\path{#1}}}
\providecommand{\Pubmed}[1]{\href{pmid:#1}{\path{#1}}}
\providecommand{\bibinfo}[2]{#2}
\ifx\xfnm\relax \def\xfnm[#1]{\unskip,\space#1}\fi
\bibitem[{Albi et~al.(2019)Albi, Bellomo, Fermo, Ha, Kim, Pareschi, Poyato and
  Soler}]{albi2019vehicular}
\bibinfo{author}{Albi, G.}, \bibinfo{author}{Bellomo, N.},
  \bibinfo{author}{Fermo, L.}, \bibinfo{author}{Ha, S.Y.},
  \bibinfo{author}{Kim, J.}, \bibinfo{author}{Pareschi, L.},
  \bibinfo{author}{Poyato, D.}, \bibinfo{author}{Soler, J.},
  \bibinfo{year}{2019}.
\newblock \bibinfo{title}{Vehicular traffic, crowds, and swarms: From kinetic
  theory and multiscale methods to applications and research perspectives}.
\newblock \bibinfo{journal}{Mathematical Models and Methods in Applied
  Sciences} \bibinfo{volume}{29}, \bibinfo{pages}{1901--2005}.
\bibitem[{Albi et~al.(2015)Albi, Herty and Pareschi}]{albi2015kinetic}
\bibinfo{author}{Albi, G.}, \bibinfo{author}{Herty, M.},
  \bibinfo{author}{Pareschi, L.}, \bibinfo{year}{2015}.
\newblock \bibinfo{title}{Kinetic description of optimal control problems and
  applications to opinion consensus}.
\newblock \bibinfo{journal}{Communications in Mathematical Sciences}
  \bibinfo{volume}{13}, \bibinfo{pages}{1407--1429}.
\bibitem[{Albi et~al.(2014)Albi, Pareschi and Zanella}]{albi2014boltzmann}
\bibinfo{author}{Albi, G.}, \bibinfo{author}{Pareschi, L.},
  \bibinfo{author}{Zanella, M.}, \bibinfo{year}{2014}.
\newblock \bibinfo{title}{Boltzmann-type control of opinion consensus through
  leaders}.
\newblock \bibinfo{journal}{Philosophical Transactions of the Royal Society A:
  Mathematical, Physical and Engineering Sciences} \bibinfo{volume}{372},
  \bibinfo{pages}{20140138}.
\bibitem[{Babovsky and Illner(1989)}]{babovsky1989convergence}
\bibinfo{author}{Babovsky, H.}, \bibinfo{author}{Illner, R.},
  \bibinfo{year}{1989}.
\newblock \bibinfo{title}{A convergence proof for {N}anbu’s simulation method
  for the full {Boltzmann} equation}.
\newblock \bibinfo{journal}{{SIAM Journal on Numerical Analysis}}
  \bibinfo{volume}{26}, \bibinfo{pages}{45--65}.
\bibitem[{Babovsky and Neunzert(1986)}]{babovsky1986simulation}
\bibinfo{author}{Babovsky, H.}, \bibinfo{author}{Neunzert, H.},
  \bibinfo{year}{1986}.
\newblock \bibinfo{title}{On a simulation scheme for the {Boltzmann} equation}.
\newblock \bibinfo{journal}{{Mathematical Methods in the Applied Sciences}}
  \bibinfo{volume}{8}, \bibinfo{pages}{223--233}.
\bibitem[{Ben~Abdallah et~al.(1996)Ben~Abdallah, Degond and
  G{\'e}nieys}]{ben1996energy}
\bibinfo{author}{Ben~Abdallah, N.}, \bibinfo{author}{Degond, P.},
  \bibinfo{author}{G{\'e}nieys, S.}, \bibinfo{year}{1996}.
\newblock \bibinfo{title}{{An energy-transport model for semiconductors derived
  from the Boltzmann equation}}.
\newblock \bibinfo{journal}{{Journal of Statistical Physics}}
  \bibinfo{volume}{84}, \bibinfo{pages}{205--231}.
\bibitem[{Bird(1970)}]{bird1970direct}
\bibinfo{author}{Bird, G.}, \bibinfo{year}{1970}.
\newblock \bibinfo{title}{Direct simulation and the {B}oltzmann equation}.
\newblock \bibinfo{journal}{The Physics of Fluids} \bibinfo{volume}{13},
  \bibinfo{pages}{2676--2681}.
\bibitem[{Burini et~al.(2016)Burini, De~Lillo and
  Gibelli}]{burini2016collective}
\bibinfo{author}{Burini, D.}, \bibinfo{author}{De~Lillo, S.},
  \bibinfo{author}{Gibelli, L.}, \bibinfo{year}{2016}.
\newblock \bibinfo{title}{Collective learning modeling based on the kinetic
  theory of active particles}.
\newblock \bibinfo{journal}{{Physics of Life Reviews}} \bibinfo{volume}{16},
  \bibinfo{pages}{123--139}.
\bibitem[{Caflisch et~al.(2021)Caflisch, Silantyev and
  Yang}]{caflisch2021adjoint}
\bibinfo{author}{Caflisch, R.}, \bibinfo{author}{Silantyev, D.},
  \bibinfo{author}{Yang, Y.}, \bibinfo{year}{2021}.
\newblock \bibinfo{title}{Adjoint {DSMC} for nonlinear {B}oltzmann equation
  constrained optimization}.
\newblock \bibinfo{journal}{Journal of Computational Physics}
  \bibinfo{volume}{439}, \bibinfo{pages}{110404}.
\bibitem[{Cordier et~al.(2005)Cordier, Pareschi and
  Toscani}]{cordier2005kinetic}
\bibinfo{author}{Cordier, S.}, \bibinfo{author}{Pareschi, L.},
  \bibinfo{author}{Toscani, G.}, \bibinfo{year}{2005}.
\newblock \bibinfo{title}{On a kinetic model for a simple market economy}.
\newblock \bibinfo{journal}{Journal of Statistical Physics}
  \bibinfo{volume}{120}, \bibinfo{pages}{253--277}.
\bibitem[{Davis et~al.(2020)Davis, Giannakis, Stadler, Stechmann and
  Manucharyan}]{davis2020super}
\bibinfo{author}{Davis, A.D.}, \bibinfo{author}{Giannakis, D.},
  \bibinfo{author}{Stadler, G.}, \bibinfo{author}{Stechmann, S.N.},
  \bibinfo{author}{Manucharyan, G.}, \bibinfo{year}{2020}.
\newblock \bibinfo{title}{Super-parameterization of {L}agrangian sea ice
  dynamics using the {B}oltzmann equation}, in: \bibinfo{booktitle}{AGU Fall
  Meeting Abstracts}, pp. \bibinfo{pages}{C048--03}.
\bibitem[{Gamba et~al.(2017)Gamba, Haack, Hauck and Hu}]{gamba2017fast}
\bibinfo{author}{Gamba, I.M.}, \bibinfo{author}{Haack, J.R.},
  \bibinfo{author}{Hauck, C.D.}, \bibinfo{author}{Hu, J.},
  \bibinfo{year}{2017}.
\newblock \bibinfo{title}{A fast spectral method for the {B}oltzmann collision
  operator with general collision kernels}.
\newblock \bibinfo{journal}{{SIAM Journal on Scientific Computing}}
  \bibinfo{volume}{39}, \bibinfo{pages}{B658--B674}.
\bibitem[{Guana et~al.(2022)Guana, Noguchia, Matsushimaa and
  Yamadaa}]{guanatopology}
\bibinfo{author}{Guana, K.}, \bibinfo{author}{Noguchia, Y.},
  \bibinfo{author}{Matsushimaa, K.}, \bibinfo{author}{Yamadaa, T.},
  \bibinfo{year}{2022}.
\newblock \bibinfo{title}{Topology optimization for rarefied gas flow problems
  using density method and adjoint {DSMC}}.
\newblock \bibinfo{journal}{Preprint.} \URLprefix
  \url{https://doi.org/10.51094/jxiv.98}.
\bibitem[{Kanazawa et~al.(2018)Kanazawa, Sueshige, Takayasu and
  Takayasu}]{kanazawa2018derivation}
\bibinfo{author}{Kanazawa, K.}, \bibinfo{author}{Sueshige, T.},
  \bibinfo{author}{Takayasu, H.}, \bibinfo{author}{Takayasu, M.},
  \bibinfo{year}{2018}.
\newblock \bibinfo{title}{{Derivation of the Boltzmann equation for financial
  Brownian motion: Direct observation of the collective motion of
  high-frequency traders}}.
\newblock \bibinfo{journal}{Physical Review Letters} \bibinfo{volume}{120},
  \bibinfo{pages}{138301}.
\bibitem[{Mohamed et~al.(2020)Mohamed, Rosca, Figurnov and
  Mnih}]{mohamed2019monte}
\bibinfo{author}{Mohamed, S.}, \bibinfo{author}{Rosca, M.},
  \bibinfo{author}{Figurnov, M.}, \bibinfo{author}{Mnih, A.},
  \bibinfo{year}{2020}.
\newblock \bibinfo{title}{Monte {C}arlo gradient estimation in machine
  learning}.
\newblock \bibinfo{journal}{Journal of Machine Learning Research}
  \bibinfo{volume}{21}, \bibinfo{pages}{1--62}.
\bibitem[{Naesseth et~al.(2017)Naesseth, Ruiz, Linderman and Blei}]{Blei2017}
\bibinfo{author}{Naesseth, C.}, \bibinfo{author}{Ruiz, F.},
  \bibinfo{author}{Linderman, S.}, \bibinfo{author}{Blei, D.},
  \bibinfo{year}{2017}.
\newblock \bibinfo{title}{Reparameterization gradients through
  acceptance-rejection sampling algorithms}.
\newblock \bibinfo{journal}{Proceedings of the 20th International Conference on
  Artificial Intelligence and Statistics (AISTATS)} .
\bibitem[{Nanbu(1980)}]{nanbu1980direct}
\bibinfo{author}{Nanbu, K.}, \bibinfo{year}{1980}.
\newblock \bibinfo{title}{Direct simulation scheme derived from the {Boltzmann}
  equation. {I.} monocomponent gases}.
\newblock \bibinfo{journal}{Journal of the Physical Society of Japan}
  \bibinfo{volume}{49}, \bibinfo{pages}{2042--2049}.
\bibitem[{Pareschi and Russo(2001)}]{pareschi2001introduction}
\bibinfo{author}{Pareschi, L.}, \bibinfo{author}{Russo, G.},
  \bibinfo{year}{2001}.
\newblock \bibinfo{title}{{An introduction to Monte Carlo method for the
  Boltzmann equation}}, in: \bibinfo{booktitle}{ESAIM: Proceedings},
  \bibinfo{organization}{EDP Sciences}. pp. \bibinfo{pages}{35--75}.
\bibitem[{Pareschi and Toscani(2013)}]{pareschi2013interacting}
\bibinfo{author}{Pareschi, L.}, \bibinfo{author}{Toscani, G.},
  \bibinfo{year}{2013}.
\newblock \bibinfo{title}{Interacting multiagent systems: kinetic equations and
  Monte Carlo methods}.
\newblock \bibinfo{publisher}{OUP Oxford}.
\bibitem[{Pareschi and Toscani(2014)}]{pareschi2014wealth}
\bibinfo{author}{Pareschi, L.}, \bibinfo{author}{Toscani, G.},
  \bibinfo{year}{2014}.
\newblock \bibinfo{title}{{Wealth distribution and collective knowledge: a
  Boltzmann approach}}.
\newblock \bibinfo{journal}{Philosophical Transactions of the Royal Society A:
  Mathematical, Physical and Engineering Sciences} \bibinfo{volume}{372},
  \bibinfo{pages}{20130396}.
\bibitem[{Poupaud(1991)}]{poupaud1991diffusion}
\bibinfo{author}{Poupaud, F.}, \bibinfo{year}{1991}.
\newblock \bibinfo{title}{{Diffusion approximation of the linear semiconductor
  Boltzmann equation: analysis of boundary layers}}.
\newblock \bibinfo{journal}{Asymptotic Analysis} \bibinfo{volume}{4},
  \bibinfo{pages}{293--317}.
\bibitem[{Villani(2002)}]{villani2002review}
\bibinfo{author}{Villani, C.}, \bibinfo{year}{2002}.
\newblock \bibinfo{title}{A review of mathematical topics in collisional
  kinetic theory}, in: \bibinfo{editor}{Friedlander, S.},
  \bibinfo{editor}{Serre, D.} (Eds.), \bibinfo{booktitle}{Handbook of
  Mathematical Fluid Dynamics}. volume~\bibinfo{volume}{1}, pp.
  \bibinfo{pages}{71--305}.
\bibitem[{Wang et~al.(2008)Wang, Lin, Caflisch, Cohen and
  Dimits}]{wang2008particle}
\bibinfo{author}{Wang, C.}, \bibinfo{author}{Lin, T.},
  \bibinfo{author}{Caflisch, R.}, \bibinfo{author}{Cohen, B.I.},
  \bibinfo{author}{Dimits, A.M.}, \bibinfo{year}{2008}.
\newblock \bibinfo{title}{Particle simulation of {C}oulomb collisions:
  Comparing the methods of {T}akizuka \& {A}be and {N}anbu}.
\newblock \bibinfo{journal}{Journal of Computational Physics}
  \bibinfo{volume}{227}, \bibinfo{pages}{4308--4329}.

\end{thebibliography}

\end{document}